\documentclass{article}

\usepackage{hyperref}
\usepackage{fullpage, calc, array, graphicx, epsfig, color, amsthm, amssymb, amsmath, mathrsfs, bigstrut}
\usepackage[all]{xy}

\input macros.tex

\newcounter{lecture}
\newtheorem{theorem}{Theorem}[lecture]
\newtheorem{proposition}[theorem]{Proposition}

\newtheorem{lemma}[theorem]{Lemma}
\newtheorem{corollary}[theorem]{Corollary}
\newtheorem{definition}[theorem]{Definition}

\newtheorem*{exercise}{Exercise}
\newtheorem*{problem}{Open Problem}
\newenvironment{examples}{\vspace{0.15cm}\noindent \textbf{Examples.}}{}
\newenvironment{examplex}{\vspace{0.15cm}\noindent \textbf{Example.}}{}
\newenvironment{example}[1]{\vspace{0.15cm}\noindent \textbf{Example: #1.}}{}
\newenvironment{remark}{\vspace{0.15cm}\noindent \textbf{Remark.}}{}

\newcommand{\NPos}{$\mathscr{N}$-position}
\newcommand{\PPos}{$\mathscr{P}$-position}
\newcommand{\mex}{\mathrm{mex}}
\newcommand{\A}{\mathscr{A}}
\newcommand{\B}{\mathscr{B}}
\newcommand{\Q}{\mathcal{Q}}
\newcommand{\T}{\mathcal{T}}
\renewcommand{\P}{\mathcal{P}}
\newcommand{\QP}{(\Q,\P)}
\newcommand{\cl}{\mathrm{cl}}
\newcommand{\hcl}{\mathrm{hcl}}
\newcommand{\sh}{\textrm{\raisebox{1pt}{\tiny \#}}}

\renewcommand{\>}{\rangle}
\newcommand{\Star}{\ast}
\mathchardef\Neg     ="0200

\newcommand{\lecture}[4]{
\stepcounter{lecture}
\begin{center}
\framebox[6.5in]{
\centering
\begin{minipage}{6.3in}
  {\bf Mis\`{e}re Games and Mis\`{e}re  Quotients} \hfill #2
  \begin{center}
    {\Large Lecture #1: #3} \\[3mm]
  \end{center}
Instructor: Aaron Siegel \hfill Scribes: #4
\end{minipage}
}
\end{center}
\bigskip
}

\begin{document}

\bibliographystyle{abbrv}

\title{Mis\`ere Games and Mis\`ere Quotients\\{\large \scshape Version 1.0}}
\author{Aaron N. Siegel}
\date{\today}
\maketitle

These notes are based on a short course offered at the Weizmann Institute of Science in Rehovot, Israel, in November 2006.  The notes include an introduction to impartial games, starting from the beginning; the basic mis\`ere quotient construction; a proof of the Guy--Smith--Plambeck Periodicity Theorem; and statements of some recent results and open problems in the subject.

First and foremost, I wish to thank the scribes for the course: Gideon Amir, Shiri Chechik, Omer Kadmiel, Amir Kantor, Dan Kushnir, Shai Lubliner, Ohad Manor, Leah Nutman, Menachem Rosenfeld, and Rivka Taub.  I also wish to thank Professor Aviezri Fraenkel for inviting me to the Weizmann Institute and suggesting this course, and thereby making these notes possible.  Finally, I wish to thank Thane Plambeck, for recognizing the importance of mis\`ere quotients and inventing this beautiful and fascinating theory.

\newpage

\section*{Introduction}

This course is concerned with \emph{impartial combinatorial games}, and in particular with mis\`ere play of such games.  Loosely speaking, a \emph{combinatorial game} is a two-player game with no hidden information and no chance elements.  We usually impose one of two winning conditions on a combinatorial game: under \emph{normal play}, the player who makes the last move wins; and under \emph{mis\`ere play}, the player who makes the last move loses.  We will shortly give more precise definitions.

The study of combinatorial games began in 1902, with C.~L.~Bouton's published solution to the game of \textsc{Nim}~\cite{bouton_1902}.  Further progress was sporadic until the 1930s, when R.~P.~Sprague~\cite{sprague_1935,sprague_1937} and P.~M.~Grundy~\cite{grundy_1939} independently generalized Bouton's result to obtain a complete theory for normal-play impartial games.

In a seminal 1956 paper~\cite{guy_1956}, R.~K.~Guy and C.~A.~B.~Smith introduced a wide class of impartial games known as \emph{octal games}, together with some general techniques for analyzing them in normal play.  Guy and Smith's techniques proved to be enormously powerful in finding normal-play solutions for such games, and they are still in active use today~\cite{flammenkamp_www_octal}.

At exactly the same time (and, in fact, in exactly the same issue of the \emph{Proceedings of the Cambridge Philosophical Society}), Grundy and Smith published a paper on mis\`ere games~\cite{grundy_1956}.  They noted that mis\`ere play appears to be quite difficult, in sharp contrast to the great success of the Guy--Smith techniques.

Despite these complications, Grundy remained optimistic that the Sprague--Grundy theory could be generalized in a meaningful way to mis\`ere play.  These hopes were dashed in the 1970s, when Conway~\cite{conway_1976} showed that the Grundy--Smith complications are intrinsic.  Conway's result shows that the most natural mis\`ere-play generalization of the Sprague--Grundy theory is hopelessly complicated, and is therefore essentially useless in all but a few simple cases.\footnote{Despite its apparent \emph{uselessness}, Conway's theory is actually quite \emph{interesting} from a theoretical point of view.  We will not say much about it in this course, but it is well worth exploring; see~\cite{conway_1976} for discussion.}

The next major advance occurred in 2004, when Thane Plambeck~\cite{plambeck_2005} recovered a tractable theory by \emph{localizing} the Sprague--Grundy theory to various restricted sets of mis\`ere games.  Such localizations are known as \emph{mis\`ere quotients}, and they will be the focus of this course.  While some of the ideas behind the quotient construction are present in Conway's work of the 1970s, it was Plambeck who recognized that the construction can be made systematic---in particular, he showed that the Guy--Smith \emph{Periodicity Theorem} can be generalized to the local setting.

This course is a complete introduction to the theory of mis\`ere quotients, starting with the basic definitions of combinatorial game theory and a proof of the Sprague--Grundy Theorem.  We include a full proof of the Guy--Smith--Plambeck Periodicity Theorem and many motivating examples.  The final lecture includes a discussion of major open problems and promising directions for future research.

\newpage
\lecture{1}{November 26, 2006}{Normal Play}{Leah Nutman \& Dan Kushnir}


\subsection*{Impartial Combinatorial Games---A Few Examples}

A combinatorial game is a two player game with no hidden information
(i.e.\ both players have full information of the game's position)
and no chance elements (given a player's move, the next position of
the game is completely determined). Let us demonstrate this notion
with a few useful examples.

\begin{figure}[ht!]
\centerline{\includegraphics[scale=0.4]{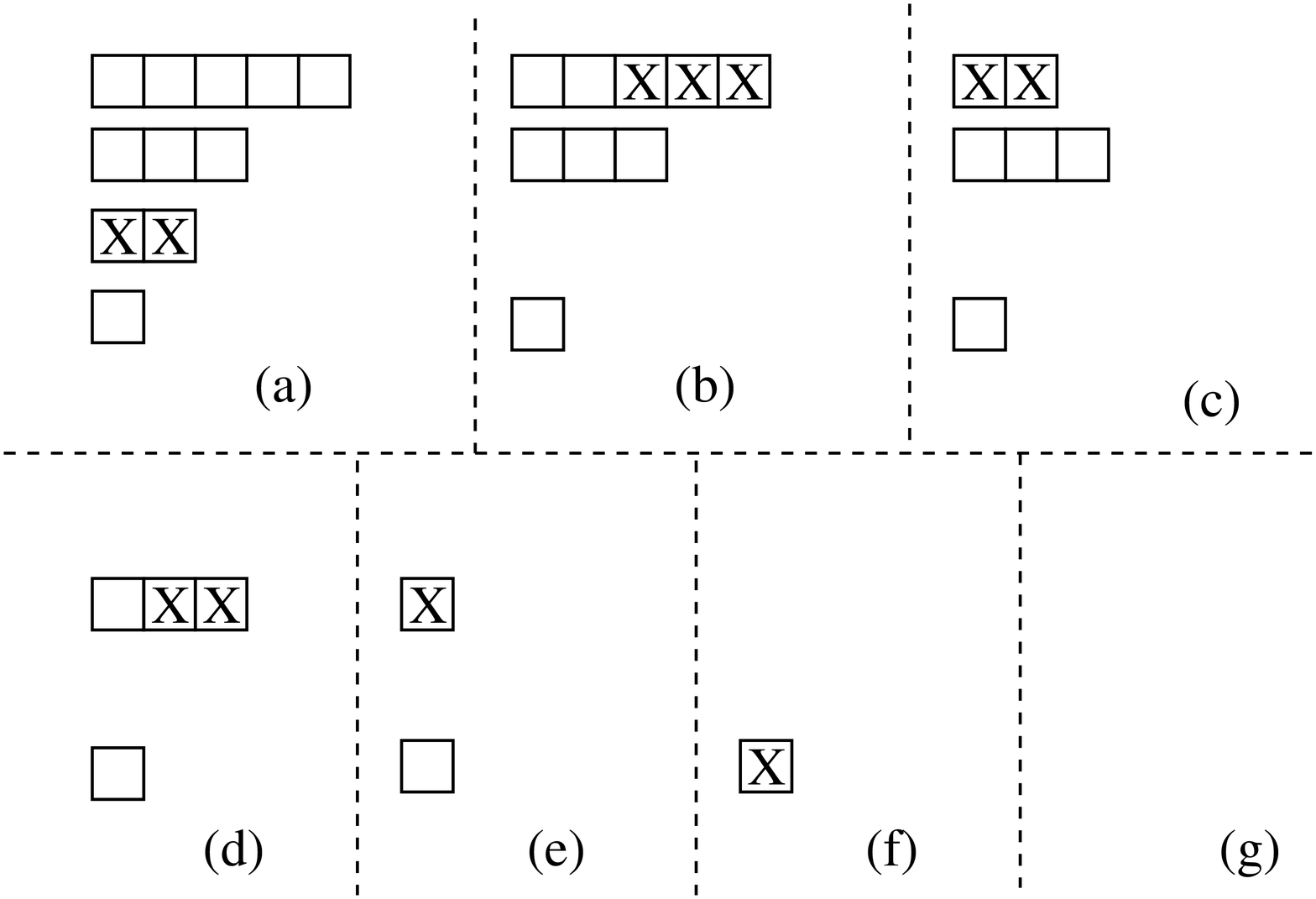}} \caption{The seven
sub-figures represent seven consecutive positions in a play of
a game of \textsc{Nim}. The last position is the empty one (with no boxes
left and thus no more possible moves). The first six mini-figures
also indicate the move taken next (which transforms the current
position into the next position): a box marked with `X' is a box
that was selected to be taken by the player whose turn it is to
play.}\label{fig:nim}
\end{figure}

\begin{example}{\textsc{\mdseries Nim}}
A position of \textsc{Nim} consists of several strips, each containing
several boxes. A move consists of removing one or more boxes from
{\em a single} strip. Whoever takes the last box (from the last
remaining strip) wins. A sample game of \textsc{Nim} is illustrated in
Figure~\ref{fig:nim}.
\end{example}

\begin{figure}[ht!]
\centerline{\includegraphics[scale=0.4]{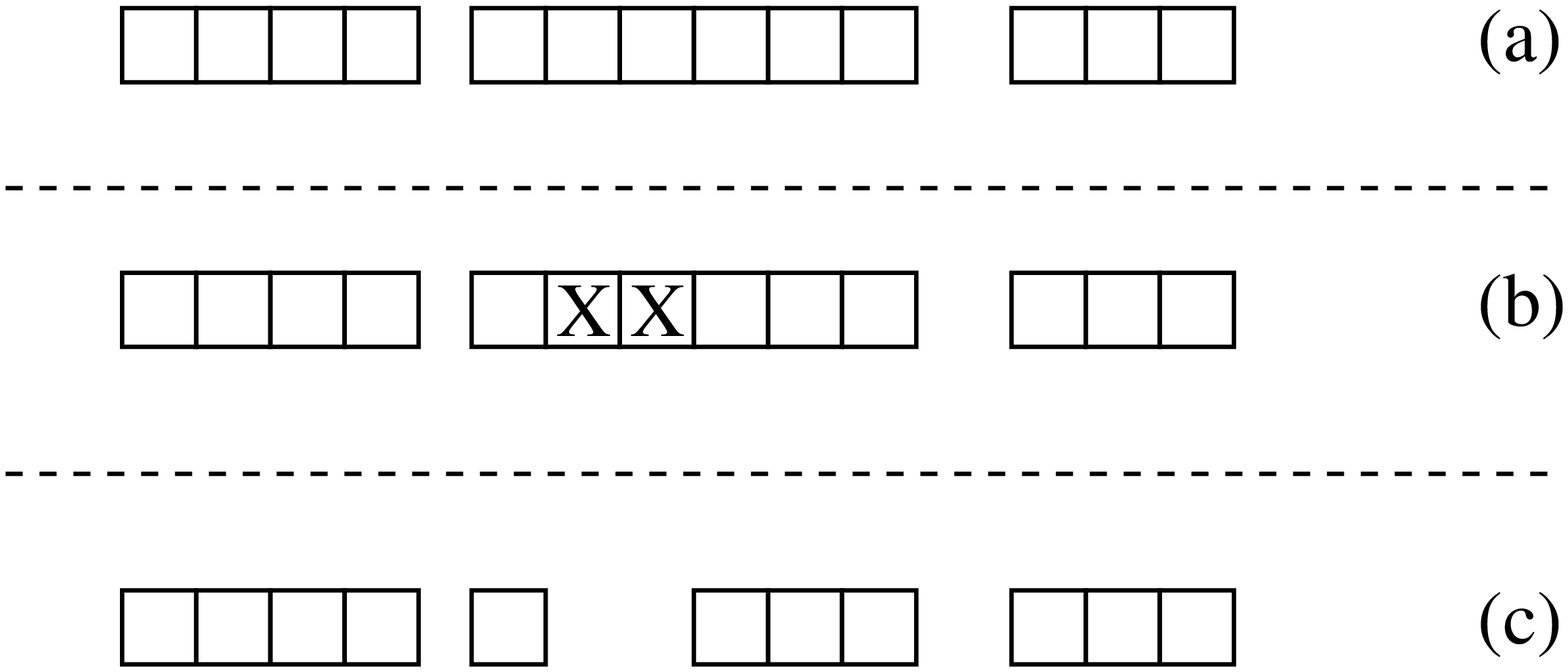}} \caption{(a)
The position before the move -- consists of three strips. (b) The
move -- the selected boxes are marked with `X'. (c) The new position
-- the middle strip was split, leaving four strips.}\label{fig:kayless}
\end{figure}

\begin{example}{\textsc{\mdseries Kayles}}
A position of \textsc{Kayles} consists of several strips, each containing
several boxes, as in \textsc{Nim}. A move consists of removing one or two \emph{adjacent} boxes
from a single strip. If the player takes a box (or two) from the
middle of a strip then this strip is split into two {\em
separate} strips. In particular, no future move can affect both
sides of the original strip. (See Figure~\ref{fig:kayless} for an
illustration of one such move.)

Whoever takes the last box (from the last remaining strip) wins.
\end{example}

\begin{example}{\textsc{\mdseries Dawson's Kayles}}
This game is identical to \textsc{Kayles} up to two differences: (1) A move
consists of removing {\em exactly two} adjacent boxes from a single
strip. (2) The winning condition is flipped: Whoever makes the last
move {\em loses}.
\end{example}

\subsection*{Winning Conditions and the Difficulty of a Game} All three
examples above share some common properties.  They are:
\begin{description}
\item[Finite.] For any given first position, there are only
finitely many possible positions that the game may take (throughout
its execution).
\item[Loopfree.] No position can occur twice in an execution of a game.
Once we leave a position, this position will never repeat itself.
\item[Impartial.] Both players have the same moves available
at all times.
\end{description}
All of the games we will consider in this course have these three
properties. As we will further discuss below, the first two
properties (finite and loopfree) imply that one of the players must
have a perfect winning strategy---that is, a strategy that guarantees
a win no matter what his opponent does.

\begin{quote}
\textbf{Main Goal:} Given a combinatorial game $\Gamma$, find an efficient winning
strategy for $\Gamma$.\footnote{More precisely, we seek a winning strategy that can be computed in polynomial time (measured against the size of a \emph{succinct description} of a game position).  In general, any use of the word ``efficient'' in this course can be safely interpreted to mean ``polynomial-time,'' though we will be intentionally vague about issues of complexity.}
\end{quote}


We will consider in this course two possible winning conditions for
our games:
\begin{description}
\item[Normal Play:] Whoever makes the last move wins.
\item[Mis\`{e}re Play:] Whoever makes the last move loses.
\end{description}
The different winning conditions of the aforementioned games turn out
to have a great effect on their difficulty. \textsc{Nim} was solved in
1902 and \textsc{Kayles} was solved in 1956.  By contrast, the solution
to \textsc{Dawson's Kayles} remains an open problem after 70 years.
(That is, we still do not know an efficient winning strategy for it.)

What makes \textsc{Dawson's Kayles} so much harder? It is exactly the fact
that the last player to move loses. In general, games with
mis\`{e}re play tend to be vastly more difficult. The themes for
this course are:
\begin{enumerate}
\item Why is mis\`{e}re play more difficult?
\item How can we tackle this difficulty?
\end{enumerate}

\subsection*{Game Representations and Outcomes}

We have mentioned that our goal is to obtain efficient winning
strategies for impartial combinatorial games. We will in fact be even more concerned
with the structure of individual positions. Therefore, by a ``game'' $G$, we
will usually mean an individual position in a combinatorial game.

Sometimes we will shamelessly abuse terminology and use the term ``game''
to refer to a system of rules.  It will (hopefully) always be clear from
the context which meaning is intended.  To help minimize confusion, we will
always denote individual positions by roman letters ($G$, $H$, $\ldots$) and
systems of rules by $\Gamma$.

\begin{figure}[ht!]
\centerline{\includegraphics[scale=0.4]{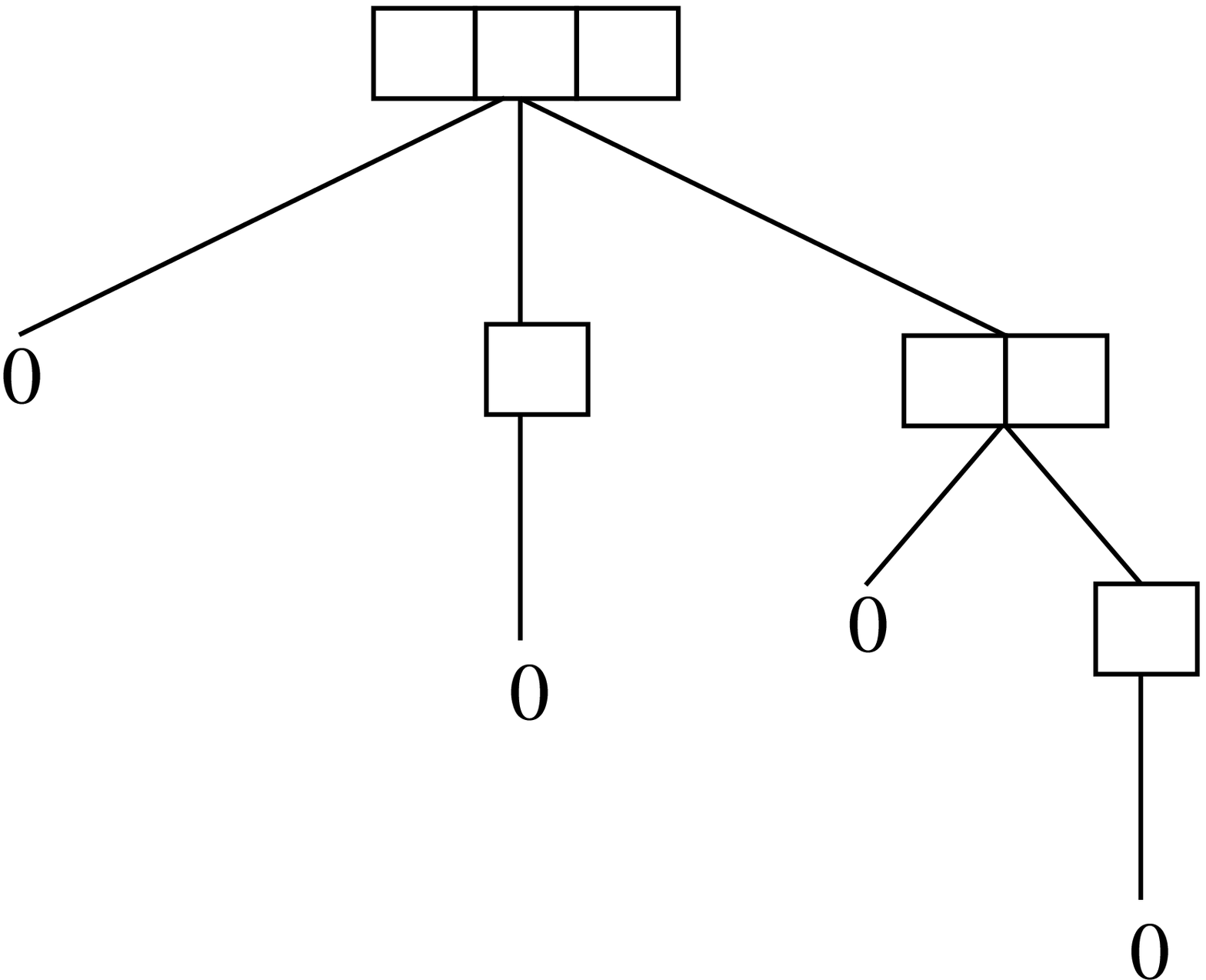}} \caption{Tree
representation of the \textsc{Nim} position $G$ which contains three boxes in
a single strip}\label{fig:NIMtree}
\end{figure}

One way to formally represent a game is as a tree. For example, the
\textsc{Nim} position $G$ which contains three boxes in a single strip can
be represented by the tree given in Figure~\ref{fig:NIMtree}.

\begin{definition}
Two games $G$ and $H$ are identical (isomorphic) if they have
isomorphic trees.  If $G$ and $H$ are identical, we write $G \cong H$.
\end{definition}

We can also think of the \textsc{Nim} position $G$ from
Figure~\ref{fig:NIMtree} as a \emph{set}: $\Box\Box\Box = \{0,\Box,\Box\Box\}$,
where $0$ denotes the game $\{\}$ with no possible moves. We call the positions we can move to directly
from a game $G$ the {\bf options} of $G$.  So we are identifying $G$ with the set of its options.

We will now introduce some notation that will make it easier to
discuss the value of any given position of a game and in particular,
the values of \textsc{Nim} positions.

\begin{definition} For every $n \geq 0$ we denote by $*n$ a strip in \textsc{Nim} of
length $n$. We write $0$ and $*$ as shorthand for $*0$ and $*1$, respectively.
Formally, we have
\[*n = \{0, * , *2, *3, \ldots, * (n-1)\}.\]
\end{definition}

As we mentioned above, every game with the properties we have
specified has a well-defined outcome (indicating who will win when
both players play perfectly). Assuming both players play perfectly,
either:
\begin{enumerate}
\item  The first player has a winning move, or
\item  Any move the first player may make will move to a position where he loses. In this case the second player can win.
\end{enumerate}

\begin{definition}
Let $G$ be a game.  The \emph{normal outcome} $o^+(G)$ is defined by
\begin{itemize}
\item $o^+(G) = \mathscr{P}$ if second player can win $G$, assuming normal play;
\item $o^+(G) = \mathscr{N}$ if first player can win $G$, assuming normal play.
\end{itemize}
Likewise, the \emph{mis\`ere outcome} $o^-(G)$ is defined by
\begin{itemize}
\item $o^-(G) = \mathscr{P}$ if second player can win $G$, assuming mis\`{e}re
play;
\item $o^-(G) = \mathscr{N}$ if first player can win $G$, assuming mis\`ere play.
\end{itemize}
We say $G$ is a \emph{normal $\mathscr{P}$-position} if $o^+(G) = \mathscr{P}$, etc.
\end{definition}

Note that $o^+$ and $o^-$ have simple recursive descriptions:
\[o^+(G) = \mathscr{P} \iff o^+(G') = \mathscr{N} \textrm{ for every option } G' \textrm{ of } G.\]
\[o^-(G) = \mathscr{P} \iff G \neq 0 \textrm{ and } o^-(G') = \mathscr{N} \textrm{ for every option } G' \textrm{ of } G.\]

$\mathscr{P}$ and $\mathscr{N}$ are short for \emph{p}revious player and \emph{n}ext player, respectively.

For example, we can consider \textsc{Nim} played with a single strip and see which positions are $\mathscr{P}$-positions and which are $\mathscr{N}$-positions:
\begin{itemize}
\item $o^+(0) = \mathscr{P}$: If there are no more boxes, then the previous move was the winning move (the previous player
took the last box).
\item $o^+(*n) = \mathscr{N}$ for every $n>0$: When there is only one strip left, the next player
can take all the remaining boxes and thus win.
\end{itemize}

What about the mis\`{e}re outcomes?
\begin{itemize}
\item $o^-(0) = \mathscr{N}$: If there are no more boxes, then the previous player
took the last box and lost. So the next player is the winning one.
\item $o^-(*) = \mathscr{P}$: When there is only one box left, the next player
must take it and lose, so the previous player is the winning one.
\item $o^-(*n) = \mathscr{N}$ for every $n>1$: Here the winning move is to take all boxes but
one.
\end{itemize}

We now revise our main goal.

\begin{quote}
\textbf{Main Goal (Revised):} Given a position $G$ in a combinatorial game, find an efficient way to compute the outcome of~$G$.
\end{quote}

In all the examples we consider in this course, the two goals are equivalent: efficient methods for computing the outcomes of positions will instantly yield efficient winning strategies.

\subsection*{Disjunctive Sums}

The positions in each of our examples naturally decompose.  In \textsc{Nim}, no single move may affect more than one strip, so each strip is effectively independent.  Both \textsc{Kayles} and \textsc{Dawson's Kayles} exhibit an even stronger form of decomposition: a typical move cuts a strip into two components, and since the components are no longer adjacent, no subsequent move can affect them both.

These observations motivate the following definition.

\begin{definition}
Let $G$ and $H$ be games.  The \emph{(disjunctive) sum} of $G$ and $H$, denoted $G+H$, is the game played as follows.  Place copies of $G$ and $H$ side-by-side.  A move consists of choosing exactly one component and making a move in that component.
\end{definition}

Formally, we can define $G + H$ as the direct sum of the trees for $G$ and $H$.  Or, thinking in terms of sets,
\[G + H = \{G' + H : G' \textrm{ is an option of } G\} \cup \{G + H' : H' \textrm{ is an option of } H\}.\]
In combinatorial game theory, it is customary to be lazy in our use of notation and write simply
\[G + H = \{G' + H, G + H'\}.\]

\subsection*{The Strategy for \textsc{Nim}}

Here is the strategy for \textsc{Nim}, assuming normal play: write the size of each strip in binary, and then do a bitwise XOR.  $G$ is a $\mathscr{P}$-position if and only if the result is identically $0$.  For example, the starting position of Figure~\ref{fig:nim}(a) has strips of sizes $5$, $3$, $2$ and $1$, so we can write
\[
\begin{array}{rr@{~$=$~}r}
       & 101 & 5 \\
\oplus & 11 & 3 \\
\oplus & 10 & 2 \\
\oplus & 1 & 1 \\ \hline
& \multicolumn{2}{l}{101} \bigstrut
\end{array}
\]
The result is nonzero, so Figure~\ref{fig:nim}(a) is an \NPos{} (in normal play).

We will shortly prove a stronger statement that implies this strategy.

\subsection*{Equivalence}

We would like to regard two games as equivalent if they behave the same way in any disjunctive sum. For now assume normal play.

\begin{definition}
\label{definition:normalequals}
We say $G$ and $H$ are \emph{equal}, and write $G = H$, iff
\[o^+(G+X) = o^+(H+X) \textrm{ for every combinatorial game } X.\]
\end{definition}

Note that if $G \cong H$, then necessarily $G = H$, but we will see that nonisomorphic games can be equal.

\begin{examplex}
$G + 0 = G$ for any game $G$.
\end{examplex}

\begin{proof}
Adding $0$ does not change the structure of $G$ at all.  (In fact, $G + 0 \cong G$.)
\end{proof}

\begin{examplex}
$G + G = 0$ for any game $G$.
\end{examplex}

\begin{proof}
We need to show that $X$ and $G+G+X$ have the same outcome, for any $X$.

First suppose $o^+(X) = \mathscr{P}$.  Second player can win $G + G + X$ as follows.  Whenever first player moves on $X$, just use the winning strategy there.  If first player ever moves on one of the copies of $G$, make the identical move on the other copy.  Second player will get the last move on $X$ because she is following the winning strategy there, and she will get the last move on $G+G$ by symmetry.

Conversely, if  $o^+(X) = \mathscr{N}$, then on $G+G+X$, just make a winning move on $X$ and proceed as before.
\end{proof}

\begin{examplex}
Here is a simple example to show how disjunctive sums can be useful for studying combinatorial games.
Consider a \textsc{Nim} position with strips of sizes $19$, $23$, $16$, $45$, $23$ and $19$.  By the previous argument, the two strips of size $19$ together equal $0$, as do the two strips of size $23$.  So this is equivalent to \textsc{Nim} with strips of sizes $16$ and $45$.

Exactly the same argument works for \textsc{Kayles} or \textsc{Dawson's Kayles}.
\end{examplex}

\begin{proposition}\ 

\begin{enumerate}
\item[(a)] $=$ is an equivalence relation.
\item[(b)] If $G = H$, then $G + K = H + K$.
\end{enumerate}
\end{proposition}

\begin{proof}
(a) is immediate, since equality of outcomes is an equivalence relation.  For (b), if $G = H$ then
\[o^+(G + X) = o^+(H + X) \textrm{ for all } X,\]
so in particular
\[o^+(G+(K+X)) = o^+(H+(K+X)) \textrm{ for all } X.\]
Disjunctive sum is associative, so $G + K = H + K$.
\end{proof}

\begin{proposition}
The following are equivalent, for games $G,H$:
\begin{enumerate}
\item[(i)] $G = H$
\item[(ii)] $o^+(G+H) = \mathscr{P}$
\end{enumerate}
\end{proposition}

\begin{proof}
(i) $\Rightarrow$ (ii): If $G = H$, then $G + G = G + H$.  But $G + G = 0$, so $o^+(G+H) = o^+(0) = \mathscr{P}$.

\vspace{0.15cm}\noindent
(ii) $\Rightarrow$ (i): By a symmetry argument (just like a previous example), $X$ and $G + H + X$ have the same outcome, for all $X$.  Therefore $G + H = 0$, so $G + H + H = H$.  But $H + H = 0$.
\end{proof}

\subsection*{The Sprague--Grundy Theorem}

\begin{theorem}[Sprague--Grundy]
For any game $G$, there is some $m$ such that $G = \Star m$.
\end{theorem}

We will in fact prove the following stronger statement.

\begin{definition}
Let $S$ be a finite set of non-negative integers.  We define the \emph{minimal excludant} of $S$, denoted $\mex(S)$, to be the least integer not in $S$.
\end{definition}

\begin{theorem}[Mex Rule]
\label{theorem:mexrule}
Suppose $G \cong \{\Star a_1,\ldots,\Star a_k\}$.  Then $G = \Star m$, where
\[m = \mex\{a_1,\ldots,a_k\}.\]
\end{theorem}

\begin{proof}
By a previous proposition, it suffices to show that $G + \Star m$ is a \PPos.  There are two cases.

\vspace{0.15cm}\noindent
\emph{Case 1}: First player moves in $G$.  This leaves the position $\Star a + \Star m$, where $\Star a$ is some option of $G$.  Since $m \not\in \{a_1,\ldots,a_k\}$, we necessarily have $a \neq m$.  If $a > m$, second player can move to $\Star m + \Star m$; if $a < m$, she can move to $\Star a + \Star a$.  In either case, she leaves a \PPos.

\vspace{0.15cm}\noindent
\emph{Case 2}: First player moves in $\Star m$.  This leaves $G + \Star a$, for some $a < m$.  Since $m$ is the \emph{minimal} excludant of $\{a_1,\ldots,a_k\}$, we must have $a = a_i$ for some $i$.  Therefore second player can move to $\Star a + \Star a$, a \PPos.
\end{proof}

The Sprague--Grundy theorem follows from one more ingredient.

\begin{exercise}
Prove the \emph{replacement lemma}: suppose $G = \{G_1,\ldots,G_k\}$ and suppose $G_1 = H$ for some $H$.  Then
\[G = \{H,G_2,\ldots,G_k\}.\]
\end{exercise}

\begin{proof}[Proof of Sprague--Grundy Theorem]
Write $G = \{G_1,\ldots,G_k\}$.  Inductively, we may assume that $G_1 = \Star a_1$, $\ldots$, $G_k = \Star a_k$.  By the replacement lemma, $G = \{\Star a_1,\ldots,\Star a_k\}$, and by the mex rule we are done.
\end{proof}
\newpage

\lecture{2}{November 27, 2006}{Octal Games and Mis\`ere Play}{Omer Kadmiel \& Shai Lubliner}

\noindent We introduce a broad class of games known as \emph{octal games}, and then give the definition of mis\`ere quotient.

\subsection*{Grundy Value}

In the previous lecture we showed:

\begin{itemize}
\item Assuming \emph{normal} play, if $G$ is any game, then $G=*m$ for
some $m$. If $G=\{*a_{1},\ldots,*a_{k}\}$ then $m=\mex\{ a_{1},\ldots,a_{k}\}$.
\item For any $G,H$, $o^{+}(G+H)=\mathscr{P}$ if and only if $G = H$.
\end{itemize}

We denote by $\mathscr{G}(G)$ the unique integer $m$ such that $G = \Star m$ in normal play.  $\mathscr{G}(G)$ is called the \emph{Grundy value} of $G$.

\subsection*{XOR and a Winning Strategy for (Normal-Play) Nim}

If $m,n$ integers then $m\oplus n$ denotes the binary XOR of $m$ and $n$.

\begin{theorem}
Let $a,b,c$ be integers.
\[o^{+}(*a+*b+*c)=\mathscr{P} \Longleftrightarrow a\oplus b\oplus c=0.\]
\end{theorem}

\begin{proof}[``Proof by Example'']
Consider the following example:

\begin{center}
\begin{tabular}{ccc}
&
11101001&
$a$\tabularnewline
$\oplus$&
01101111&
$b$\tabularnewline
$\oplus$&
00000111&
$c$\tabularnewline
\hline
&
10000001&
\tabularnewline
\end{tabular}
\end{center}

As the XOR of these values $\neq0$, we must show that this is an \NPos.
The first player simply finds the most significant bit marked $1$ in the XOR and chooses
any component in which this bit is a $1$.  In this example, that component is $a$.
He then makes an appropriate move in $a$ that switches the most significant bit to $0$,
and sets all lower-order bits as needed to make the sum equal $0$.  Here the winning
move is from $a$ to $a' = 01101000$, changing just the first and last bits.
$a' \oplus b \oplus c = 0$, so by induction it is a \PPos.
\end{proof}

\begin{corollary}
\label{corollary:nimaddition}
$*a+*b=*(a\oplus b)$
\end{corollary}

\begin{proof}
$a\oplus b\oplus(a\oplus b)=0$, so $*a+*b+*(a\oplus b)=0$.
\end{proof}

\subsection*{Example: \textsc{Dawson's Kayles}}

Recall that in \textsc{Dawson's Kayles}, a move consists of removing \emph{exactly}
two adjacent boxes.  We defined \textsc{Dawson's Kayles} as a mis\`ere-play game, but we can just as easily consider it in normal play.
Denote by $H_{n}$ a single strip of Length $n$.  Then the moves from~$H_n$ are to $H_a + H_{n-2-a}$, where $1 \leq a \leq n-2$.

We can use the Sprague--Grundy theorem and the \textsc{Nim} addition rule to compute normal-play values of~$H_n$ easily.

\begin{center}
\begin{tabular}{l}
$H_{0}=\{\} = 0$\\

$H_{1}=\{\} = 0$\\

$H_{2}=\{H_0\}=\{0\}=*$\\

$H_{3}=\{ H_{1}\}=\{0\}=*$\\

$H_{4}=\{ H_{1}+H_{1},H_{2}+H_{0}\}=\{0+0,*+0\}=\{0,*\}=*2$\\

$H_5 = \{H_2 + H_1, H_3 + H_0\} = \{*+0,*+0\} = \{*,*\}=0$\\

$H_6 = \{H_2+H_2,H_3+H_1,H_4+H_0\} = \{*+*,*+0,*2+0\} = \{0,*,*2\}=*3$
\end{tabular}
\end{center}

This rapidly becomes tedious, and it's easily implemented on a computer.  The results of a computer calculation are shown in Figure~\ref{figure:normaldk}.  Each row represents a block of 34 Grundy values: the first row shows $\mathscr{G}(H_0)$ through $\mathscr{G}(H_{33})$; the next row shows $\mathscr{G}(H_{34})$ through $\mathscr{G}(H_{67})$; etc.  The number 34 was obviously not chosen by accident; after a few initial anomalies, a strong regularity quickly emerges with period 34.  We now prove a theorem that shows, for a wide class of games, that if such periodicity is observed for ``sufficiently long'' (in a sense to be made precise) then it must continue forever.

\begin{figure}[h]
\centering
\begin{tabular}{r|l}
& $0{\color{white}~0~0~0~0~0~0~0~0~0~}1{\color{white}~0~0~0~0~0~0~0~0~0~}2{\color{white}~0~0~0~0~0~0~0~0~0~}3$ \\
& $0~1~2~3~4~5~6~7~8~9~0~1~2~3~4~5~6~7~8~9~0~1~2~3~4~5~6~7~8~9~0~1~2~3$\\
\hline
{$0+$} & {$0~0~1~1~2~0~3~1~1~0~3~3~2~2~4~0~5~2~2~3~3~0~1~1~3~0~2~1~1~0~4~5~2~7$}  \\
{$34+$} & {$4~0~1~1~2~0~3~1~1~0~3~3~2~2~4~4~5~5~2~3~3~0~1~1~3~0~2~1~1~0~4~5~3~7$}  \\
{$68+$} & {$4~8~1~1~2~0~3~1~1~0~3~3~2~2~4~4~5~5~9~3~3~0~1~1~3~0~2~1~1~0~4~5~3~7$}  \\
{$102+$} & {$4~8~1~1~2~0~3~1~1~0~3~3~2~2~4~4~5~5~9~3~3~0~1~1~3~0~2~1~1~0~4~5~3~7$}  \\
{$136+$} & {$4~8~1~1~2~0~3~1~1~0~3~3~2~2~4~4~5~5~9~3~3~0~1~1~3~0~2~1~1~0~4~5~3~7$}  \\
{$170+$} & {$4~8~1~1~2~0~3~1~1~0~3~3~2~2~4~4~5~5~9~3~3~0~1~1~3~0~2~1~1~0~4~5~3~7$}  \\
\end{tabular}
\caption{Grundy values of \textsc{Dawson's Kayles} in normal play.}
\label{figure:normaldk}
\end{figure}

\subsection*{Octal Games and Octal Codes}

\begin{definition}
An \emph{octal code} is a sequence of digits $0.d_{1}d_{2}d_{3}\ldots$ where $0 \leq d_i < 8$ for all $i$.
\end{definition}

An octal code specifies the rules for a particular \emph{octal game}.  An octal game is played with strips of boxes, and the code describes how many boxes may be removed and under what circumstances.  The digit~$d_k$ specifies the conditions under which $k$ boxes may be removed.

Let us consider the bit
representation of each $d_{k}$: Denote $d_{k}=\epsilon_{0}+\epsilon_1\cdot2+\epsilon_2\cdot4$, where each $\epsilon_i = 0$ or $1$.

\begin{itemize}
\item We can remove an entire strip of length $k$ iff $\epsilon_0 = 1$.
\item We can remove $k$ boxes from the \emph{end} of a strip (leaving at least one box) iff $\epsilon_1 = 1$.
\item We can remove $k$ boxes from the \emph{middle} of a strip (leaving at least one box on each end) iff $\epsilon_2 = 1$.
\end{itemize}

Therefore: 

\begin{itemize}
\item \textsc{Dawson's Kayles} is represented by $0.07$ as you have to remove
exactly two blocks every time from anywhere in the strip, and you
can remove an entire strip of length 2.
\item \textsc{Kayles} is represented by $0.77$ as you can remove one or two
boxes from a single strip.
\item \textsc{Nim} is represented by the infinite sequence $0.3333333\ldots$
as you are allowed to take any number of boxes from the end or to take an entire strip of any length (but you are not allowed to separate the original strip into two strips).
\end{itemize}

\subsection*{Guy--Smith Periodicity Theorem}

\begin{theorem}[Guy--Smith Periodicity Theorem]
Consider an octal game with \emph{finitely many} non-zero
code digits, and let $k$ be largest with $d_k \neq 0$. Denote by $H_n$ a strip of length $n$. Suppose that for some $n_{0} > 0$ and $p > 0$ we have
\[\mathscr{G}(H_{n+p})=\mathscr{G}(H_{n}) \textrm{ for every $n$ with } n_{0}\leq n<2n_{0}+p+k.\]
Then
\[\mathscr{G}(H_{n+p})=\mathscr{G}(H_{n}) \textrm{ for all } n \geq n_0.\]
\end{theorem}

\begin{proof}
Note that a move from $H_{n}$ is always to $H_a + H_b$, where $n-k \leq a+b < n$.  (In taking a whole strip, or from the end of the strip, we may take one or both of $a,b$ to be $0$.)

Now proceed by induction on $n$.  The base case $n < 2n_0+p+k$ is given by hypothesis, so assume $n \geq 2n_0+p+k$.  A move from $H_{n+p}$ is to $H_{a}+H_{b}$ where $a+b \geq n+p-k$.

Since $n \geq 2n_0+p+k$, we have $n+p-k\geq2n_{0}+2p$, so without loss of generality $b\geq n_{0}+p$.
(Since the sum $a+b$ is greater than or equal to $2(n_{0}+p)$, at least one of
the elements must be at least $n_0+p$.)
By induction $\mathscr{G}(H_{b-p})=\mathscr{G}(H_{b})$, so
\[\mathscr{G}(H_{a}+H_{b-p}) = \mathscr{G}(H_a) \oplus \mathscr{G}(H_{b-p})
= \mathscr{G}(H_a) \oplus \mathscr{G}(H_b) = \mathscr{G}(H_{a}+H_{b}).\]
Here is the picture:

\bigskip

\begin{tabular}[c]{|c|c|c|crl}
\cline{1-3}
\makebox[3cm]{}&$\times$&\makebox[4cm]{}&\makebox[1.2cm]{}&\makebox[0.5cm][l]{$H_n$}&$\to H_a + H_{b-p}$\\
\cline{1-3}
\end{tabular}

\bigskip

\begin{tabular}[c]{|c|c|c|crl}
\cline{1-3}
\makebox[3cm]{}&$\times$&\makebox[5cm]{}&\makebox[0.2cm]{}&\makebox[0.5cm][l]{$H_{n+p}$}&$\to H_a + H_b$\\
\cline{1-3}
\end{tabular}

\begin{tabular}[t]{c c c r}
$\underbrace{\makebox[3.1cm]{}}$&\makebox[0.06cm]{}&$\underbrace{\makebox[5.3cm]{}}$&\\
$a$&&$b$&\\
\end{tabular}

\bigskip

Now $H_a+H_{b-p}$
is an option of $H_{n}$, so we conclude that the options of $H_{n+p}$ have exactly
the same $\mathscr{G}$-values as those of~$H_{n}$. Since the $\mathscr{G}$-values of $H_{n+p}$ and $H_n$ both observe the mex rule, we have
\[\mathscr{G}(H_{n+p})=\mathscr{G}(H_{n}). \qedhere\]
\end{proof}

When $p$ and $n_0$ are as small as possible, we say that $\Gamma$ has (normal-play) \emph{period} $p$ and \emph{preperiod} $n_0$.

\begin{examples}
In normal play:
\begin{itemize}
\item \textsc{Kayles} (0.77) has period 12.
\item \textsc{Dawson's Kayles} (0.07) has period 34.
\item 0.106 has period 328226140474.  (See~\cite{flammenkamp_www_octal}.)
\item 0.007 is not known to be periodic.
\end{itemize}
\end{examples}

\begin{problem}
Does there exist a finite octal code (i.e., an octal code with finitely many non-zero digits) that yields an aperiodic game?
\end{problem}

\subsection*{Mis\`ere Nim}

We now consider \textsc{Nim} in mis\`ere play.  It is not hard to show the following.  If $G$ consists of heaps of sizes $a_1,\ldots,a_k$, then
\[o^{-}(G)=\mathscr{P}\Longleftrightarrow a_{1}\oplus a_{2}\oplus\cdots\oplus a_{k}=0,\]
\emph{unless} every $a_{i}$ is equal to 0 or 1.  In that case, $o^{-}(G)=\mathscr{P}\Longleftrightarrow a_{1}\oplus\cdots\oplus a_{k}=1$.

So the strategy for mis\`ere \textsc{Nim} is: play exactly like in normal \textsc{Nim}, unless
your move would leave only heaps of size 0 or 1. In that case, play
to leave an odd number of heaps of size 1.

\subsection*{Mis\`ere Equality}

We now make the exact same definition of equality as before (cf.~Definition~\ref{definition:normalequals}), this time assuming mis\`ere play.

\begin{definition}
\label{definition:equals}
\[G=H\Longleftrightarrow o^{-}(G+X)=o^{-}(H+X) \textrm{ for all } X\]
\end{definition}

Recall that in normal play any two \PPos{}s are equal (and in particular, any \PPos{} is equal to $0$). We shall
see that this is not the case in mis\`ere play.

In mis\`ere play:

\begin{itemize}
\item $0$ is an \NPos{}.
\item $*$ is a \PPos{}.
\item $*2$ is an \NPos.
\end{itemize}
This we have already seen.  Note that $*2+*2$ is also a \PPos.  No matter what first player does, second player can always respond by moving to $*$:
\[\xymatrix{
\Star 2 + \Star 2 \ar[dr] \ar[drr] \\
& \Star 2 + 0 \ar[dl] & \Star 2 + \Star \ar[dll] \\
\Star + 0 
}\]
This immediately shows that $*2+*2 \neq 0$, since $*2+*2$ is a \PPos{} but $0$ is an \NPos.  In fact, we will now show that $\Star 2 + \Star 2 \neq \Star$, thus exhibiting two distinct \PPos{}s.

\begin{proposition}
$*+*=0$.
\end{proposition}

\begin{proof}
Whoever can win $X$ can also win $X+*+*$: he follows the winning strategy on $X$, and if his opponent ever moves on one copy of $\Star$, he responds by moving on the other.  This guarantees that his opponent will make the last move on $X$, leaving either $0$ or $\Star + \Star$. But both of these are \NPos{}s.
\end{proof}

Now $\Star + \Star 2 + \Star 2$ is an \NPos, since it has a move to $\Star 2 + \Star 2$.  The following proposition therefore shows that $\Star \neq \Star 2 + \Star 2$.

\begin{proposition}
\label{proposition:fourstar2}
$*2+*2+*2+*2$ is a \PPos.
\end{proposition}

\begin{proof}
The options are $\Star2+\Star2+\Star2+0$ and $\Star2+\Star2+\Star2+\Star$.  But these have moves to $\Star2+\Star2+0+0$ and $\Star2+\Star2+\Star+\Star$, respectively.  By the previous proposition, both of these are equal to $\Star2+\Star2$, a \PPos.
\end{proof}

In fact, it is possible to show that $\Star 2 + \Star 2 \neq \Star m$ for any $m$.  So even among sums of \textsc{Nim}-heaps, we have games that are not equivalent to any \textsc{Nim}-heap.  This contrasts sharply with the situation in normal play, where \emph{every} game is equivalent to a \textsc{Nim}-heap.

We have seen that $\Star+\Star = 0$.  There are very few other identities we can establish in mis\`ere play.  Here are really the only two:

\begin{exercise}[Mis\`ere Mex Rule]
Suppose $G \cong \{\Star a_1,\ldots,\Star a_k\}$.  Then $G = \Star m$, where
\[m = \mex\{a_1,\ldots,a_k\},\]
\emph{provided that} at least one $a_i = 0$ or $1$.  (cf.~Theorem~\ref{theorem:mexrule})
\end{exercise}

\begin{exercise}
For any $m$, we have $\Star m + \Star = \Star(m \oplus 1)$.  (cf.~Corollary~\ref{corollary:nimaddition})
\end{exercise}

The mis\`ere mex rule is spectacularly false if every $a_i \geq 2$.  For example, let
\[G = \{\Star 2\},\]
the game whose only option is $\Star 2$.  ($G$ is sometimes called $\Star 2_\sh$.)  $G$ is a \PPos, so right away we have $G \neq 0$.  As an exercise, show that $G$ is not equal to any $\Star m$.  In fact, it is possible to show that $G$ is not equal to any \emph{sum} of \textsc{Nim}-heaps, but we won't do that in this course.

\subsection*{Birthdays}

Clearly, things are more complicated in mis\`ere play than in normal play.  We now state a result that shows just how much worse they are.

\begin{definition}
The \emph{birthday} of a game $G$ is the height of its game tree.
\end{definition}

In \emph{normal} play there are just six games with birthday $\leq 5$: $0$, $*$, $*2$, $*3$, $*4$, and $*5$.  In mis\`ere play, there are 4171780.  On day 6 there are more than $2^{4171779}$ \dots

The theory of mis\`ere games modulo $=$ is beautiful and fascinating, but these results suggest that it is not terribly useful: we very quickly run into seemingly intractable complications.  We will not say much more about this ``global theory'' in this course; the interested reader is referred to~\cite{conway_1976}.

\subsection*{Mis\`ere Quotients}

If $G = H$, then $G + X$ and $H + X$ have the same outcomes, for any game~$X$.  As we've just observed, this equality relation gives rise to a virtually intractable theory.  The problem is that $G = H$ is too strong a relation---we are requiring that $G$ and $H$ behave identically in any context, which is asking a bit too much.

\begin{quote}
\emph{Key Idea}: Suppose we just want to know how to play \textsc{Kayles} (for example).  We just need to specify how a \textsc{Kayles} position~$G$ interacts with other positions that actually occur in \textsc{Kayles}.
\end{quote}

With this in mind, fix a set $\A$ of games (usually, $\A$ will be the set of positions that occur in some octal game).  Assume that $\A$ is closed under addition.

\begin{definition}
Let $\A$ be a set of games, closed under addition.  Then for $G,H \in \A$,
\[G \equiv_\A H \Longleftrightarrow o^-(G+X) = o^-(H+X) \textrm{ for all } X \in \A.\]
\end{definition}

Compare this to Definition~\ref{definition:equals}: we are restricting the domain of games that can be used to distinguish $G$ from $H$.  This coarsens the equivalence and allows us to recover a tractable theory.  Very often, the set of equivalence classes modulo $\equiv_\A$ is finite, even when $\A$ is infinite.  (It is trivial to see that $\equiv_\A$ is an equivalence relation, since outcome-equality is an equivalence relation.)

Now, think of normal-play Grundy values as elements of the group
\[\mathcal{D} = \bigoplus_{\mathbb{N}} \mathbb{Z}_2,\]
a (countably) infinite direct sum of copies of $\mathbb{Z}_2$ (one for each binary digit).  The Sprague--Grundy theory maps each game~$G$ to an element of~$\mathcal{D}$, thus representing the normal-play structure of~$G$ in terms of the group structure of~$\mathcal{D}$.  We will show that the equivalence classes modulo $\equiv_\A$ function as a \emph{localized} mis\`ere analogue of the Sprague--Grundy theory.

We will make a slightly stronger assumption on $\A$ than closure under addition.

\begin{definition}
A set of games $\A$ is \emph{hereditarily closed} if, for any $G \in \A$ and any option $G'$ of $G$, we also have $G' \in \A$.
\end{definition}

\begin{definition}
$\A$ is \emph{closed} if it is both hereditarily closed and closed under addition.
\end{definition}

Note that if $\A$ is the set of positions that occur in an octal game, then $\A$ is closed.  In fact, virtually all sets of games that are interesting to us are closed, so there is little harm in making this assumption.

\begin{examplex}
Let $\A = \{$all sums of $*$ and $*2\}$, that is,
\[\A = \{m \cdot \Star + n \cdot \Star 2 : m,n \in \mathbb{N}\}.\]
Let's compute the equivalence classes modulo $\equiv_\A$.
\begin{itemize}
\item $\Star \not\equiv_\A 0$, since $\Star$ is a \PPos{} and $0$ is an \NPos.
\item Likewise, $\Star 2 \not\equiv_\A \Star$ since $\Star 2$ is an \NPos.  Further, $\Star 2 \not\equiv_\A 0$: let $X = \Star 2$; then $\Star 2 + X = \Star 2 + \Star 2$ is a \PPos, but $0 + X = \Star 2$ is an \NPos.
\item Finally, $\Star 2 + \Star \not\equiv_\A \Star$ since it's an \NPos; $\Star 2 + \Star \not\equiv_\A 0$, since they're distinguished by $X = \Star$; and $\Star 2 + \Star \not\equiv_\A \Star 2$, since they're distinguished by $X = \Star 2$.
\end{itemize}
This gives four equivalence classes:
\[
\begin{array}{cccc}
[0] & [\Star] & [\Star 2] & [\Star 2 + \Star] \bigstrut \\
\mathscr{N} & \mathscr{P} & \mathscr{N} & \mathscr{N} \bigstrut
\end{array}
\]
Are there others?
\begin{itemize}
\item Yes!  $\Star 2 + \Star 2$ is a \PPos, so it's either equivalent to $\Star$, or a new equivalence class.  But:
\begin{itemize}
\item $\Star + (\Star 2 + \Star 2)$ is an \NPos, since it has a move to $\Star 2 + \Star 2$, which is $\mathscr{P}$;
\item $\Star 2 + \Star 2 + (\Star 2 + \Star 2)$ is a \PPos{} (Proposition~\ref{proposition:fourstar2}).
\end{itemize}
Therefore $\Star \not\equiv_\A \Star 2 + \Star 2$.
\item Similar reasoning shows that $\Star 2 + \Star 2 + \Star$ gives yet another equivalence class.
\end{itemize}
So we have six equivalence classes total:
\[
\begin{array}{cccccc}
[0] & [\Star] & [\Star 2] & [\Star 2 + \Star] & [\Star 2 + \Star 2] & [\Star 2 + \Star 2 + \Star] \bigstrut \\
\mathscr{N} & \mathscr{P} & \mathscr{N} & \mathscr{N} & \mathscr{P} & \mathscr{N} \bigstrut
\end{array}
\]
We now show that these are the only six.
\end{examplex}

\begin{lemma}
Let $n \geq 1$.  Then $n \cdot \Star 2$ is a \PPos{} iff $n$ is even.
\end{lemma}

\begin{proof}
If $n$ is even, then second player's strategy is to cancel out copies of $\Star 2$ (using the fact that $\Star + \Star = 0$) until we get down to $\Star 2 + \Star 2$, which is known to be a \PPos.

If $n$ is odd, $n \geq 3$, then first player can win by moving to $(n-1) \cdot \Star 2$.

Finally, if $n = 1$, then first player simply moves to $\Star$.
\end{proof}

\begin{lemma}
Let $n \geq 1$.  Then $n \cdot \Star 2 + \Star$ is always an \NPos.
\end{lemma}

\begin{proof}
If $n$ is even, then the winning move is to $n \cdot \Star 2$, which is a \PPos{} by the previous Lemma.

If $n$ is odd, $n \geq 3$, then the winning move is to $(n-1) \cdot \Star 2 + \Star + \Star$, which again is a \PPos, since $\Star + \Star = 0$.

Finally, if $n = 1$, then the winning move is to $0 + \Star$.
\end{proof}

\begin{corollary}
Suppose $G = m \cdot \Star + n \cdot \Star 2$ and $X = m' \cdot \Star + n' \cdot \Star 2$.  If $n \geq 1$, then the outcome of $G + X$ depends only on the parities of $m+m'$ and $n+n'$.
\end{corollary}

\begin{proof}
Follows immediately from the previous two Lemmas and the fact that $\Star + \Star = 0$.
\end{proof}

\begin{corollary}
Let $G = m \cdot \Star + n \cdot \Star 2$ and $H = m' \cdot \Star + n' \cdot \Star 2$.
\[\textrm{If } n,n' \geq 1,\ m \equiv m' \ (\textrm{mod}~2), \textrm{ and } n \equiv n' \ (\textrm{mod}~2), \textrm{ then } G \equiv_\A H.\]
\end{corollary}

\begin{corollary}
There are exactly six equivalence classes modulo $\equiv_\A$.
\end{corollary}

\begin{proof}
By the previous corollary, every $G \in \A$ is equivalent to $m \cdot \Star + n \cdot \Star 2$, for some $m < 2$ and $n < 3$.  There are only six such possibilities, and we've already shown that all six are mutually inequivalent.
\end{proof}

\textbf{Warning.} We've just shown that $\Star 2+\Star 2+\Star 2 \equiv_\A \Star 2$.  However, equality does \emph{not} hold:

\begin{exercise}
Show that $\Star 2 + \Star 2 + \Star 2 \neq \Star 2$.  (Hint: try $X = \Star2_\sh 1$, defined by $\Star2_\sh1 = \{\Star2_\sh,\Star\} = \{\{\Star 2\}, \Star\}$.)
\end{exercise}

This shows that the equivalence $\equiv_\A$ is a genuine coarsening of equality.  There exist unequal games that are equivalent modulo~$\A$.

This finishes our example.  We now return to the general context.

\begin{lemma}
Let $\A$ be any closed set of games and $G,H \in \A$.  If $G \equiv_\A H$ and $K \in \A$, then $G + K \equiv_\A H + K$.
\end{lemma}

\begin{proof}
For $X \in \A$, we have
\[o^-((G+K) + X) = o^-(G+(K+X)) \textrm{ and } o^-(H+(K+X)) = o^-((H+K)+X).\]
But $\A$ is closed, so $K + X \in \A$.  Since $G \equiv_\A H$, we have $o^-(G+(K+X)) = o^-(H+(K+X))$, as needed.
\end{proof}

Moreover, since $\A$ is hereditarily closed, we have $0 \in \A$.  So the equivalence class of $0$ is an identity, and in fact we have a monoid.

\begin{definition}
A \emph{semigroup} is a set $S$ equipped with an associative binary operation $\cdot$.  That is,
\begin{itemize}
\item If $x,y \in S$, then $x \cdot y \in S$;
\item If $x,y,z \in S$, then $x \cdot (y \cdot z) = (x \cdot y) \cdot z$.
\end{itemize}
A semigroup $S$ is a \emph{monoid} if it has an identity, and \emph{commutative} if its operation is commutative.
\end{definition}

We've shown that the equivalence classes of $\A$ modulo $\equiv_\A$ form a commutative monoid $\Q$.
\[\Q = \{[G]_{\equiv_\A} : G \in \A\}.\]
Furthermore, if $G \equiv_\A H$, then since $o^-(G+0) = o^-(H+0)$, we have
\[G \textrm{ is a \PPos} \Longleftrightarrow H \textrm{ is a \PPos}.\]
So we can define a subset $\P \subset \Q$ by
\[\P = \{[G]_{\equiv_\A} : G \in \A \textrm{ is a \PPos}\}.\]

\begin{definition}
The structure $\QP$ is the \emph{mis\`ere quotient} of $\A$, and we denote it by $\Q(\A)$.
\end{definition}

\begin{examplex}
Let's sketch the structure of $\Q(\A)$ for our example
\[\A = \{\textrm{sums of } \Star \textrm{ and } \Star 2\}.\]
Denote by $\Phi : \A \to \Q$ the quotient map
\[\Phi(G) = [G]_{\equiv_\A}.\]
Now $\A$ is generated (as a monoid) by $\Star$ and $\Star 2$.  Put
\[1 = \Phi(0) = [0] \qquad a = \Phi(\Star) = [\Star] \qquad b = \Phi(\Star 2) = [\Star 2].\]
We know that $\Star + \Star = 0$, so in fact $a^2 = 1$.  Furthermore, we've seen that $\Star 2 + \Star 2 + \Star 2 \equiv_\A \Star 2$, so we have $b^3 = b$.  But we also know that the six elements
\[
\begin{array}{ccccccc}
\A & [0] & [\Star] & [\Star 2] & [\Star 2 + \Star] & [\Star 2 + \Star 2] & [\Star 2 + \Star 2 + \Star] \\
\downarrow \\
\Q & 1 & a & b & ab & b^2 & ab^2
\end{array}
\]
are all distinct.  Thus $\Q = \{1,a,b,ab,b^2,ab^2\}$ and we have the presentation
\[\Q \cong \<a,b~|~a^2=1,\ b^3=b\>.\]
Since $\Star$ and $\Star 2 + \Star 2$ are the only \PPos{}s (up to equivalence), we also have $\P = \{a,b^2\}$.  This mis\`ere quotient is called $\mathcal{T}_2$, and it is the first of many that we will see.
\end{examplex}

\newpage
\lecture{3}{November 28, 2006}{The Periodicity Theorem}{Amir Kantor \& Gideon Amir}

\begin{definition}[Definitions]
Let~$\A$ be any set of games. Define
\[
\hcl(\A)\triangleq\left\{ \textrm{subpositions of all games in
}\A\right\},
\]
\[
\cl(\A)\triangleq\textrm{Closure under addition of }\hcl(\A).
\]
\end{definition}

\begin{remark}
$\cl(\A)$ is hereditarily closed. To see this, let~$G=G_{1}+G_{2}+\ldots+G_{k}$ where~$G_{i}\in\hcl(\A)$.
W.l.o.g.~$G'=G_{1}'+G_{2}+\ldots+G_{k}$. We know that~$G_{1}'\in\hcl(\A)$ since the latter is hereditarily closed.
\end{remark}

\begin{examplex}
$\cl(\{\Star 2\}) = \{$sums of $\Star,\Star2\} = \{i \cdot \Star + j \cdot \Star2 : i,j \in \mathbb{N}\}$.
\end{examplex}

\begin{exercise}
\item
\begin{itemize}
\item If $\A\subseteq \B$ and~$\B$ is closed, then $\cl(\A)\subseteq \B$.
\item $\textrm{cl}(\textrm{cl}(\A))=\cl(\A)$.
\end{itemize}
\end{exercise}

\begin{definition}
If~$\A$ is not closed, $\Q(\A)\triangleq \Q(\cl(\A))$. We sometimes write
$\Q(G)\triangleq \Q\left(\cl\left(\{G\}\right)\right)$.
\end{definition}

\begin{examplex}
$\mathcal{T}_{2}\cong{}\Q(*2)$.
\end{examplex}

\subsection*{Quotients of Octal Games}

Let's consider the context of a specific octal game, such as \textsc{Kayles}.  Denote by $H_n$ a \textsc{Kayles} heap of size~$n$ and let $\A$ be the set of all \textsc{Kayles} positions; that is,
\[\A = \cl(H_0,H_1,H_2,H_3,\ldots).\]
Let $\QP$ be the mis\`ere quotient for \textsc{Kayles} and consider the quotient map $\Phi : \A \to \Q$.

\begin{remark}
If we know~$\Phi(H_n)$ for all $n$, then
if~$G=H_{n_1}+\cdots+H_{n_k}$ we can easily compute
$\Phi(G)=\Phi(H_{n_1})\cdots\Phi(H_{n_k})$.  So, in order to specify $\Phi$, it suffices to specify the single-heap values $\Phi(H_n)$.
\end{remark}

\vspace{0.15cm}\noindent
The \emph{main point} is:
\begin{quote}
Suppose we know~$\Q(\A)$, together with~$\Phi(H_n)$ for all $n$. If we
want to know~$o^{-}(G)$ for~$G\in \A$, we can
write~$G=H_{n_1}+\cdots+H_{n_k}$,
compute~$\Phi(G)=\Phi(H_{n_1})\cdots\Phi(H_{n_k})$, and simply look up whether
$\Phi(G) \in \P$. If the quotient
is finite, we've reduced the problem of finding~$o^{-}(G)$ to a
small number of operations on a finite multiplication table.  This
yields an efficient way to compute $o^-(G)$.
\end{quote}

\vspace{0.15cm}
So we direct our energies at computing the values of $\Phi(H_n)$ for all $n$.  In practice, we can construct good algorithms for computing quotients of a finite number of heaps.  (We won't have time to discuss them in this course; see~\cite[Appendix A]{siegel_200Xd}.)  If we run these algorithms on \textsc{Kayles} to heap 120, we get the result shown in Figure~\ref{figure:kayles}.

\begin{figure}
\centering
\[\begin{array}{@{}r@{~}l@{}}
\Q(H_0,H_1,H_2,\ldots,H_{120}) \cong &
\hangingpresent[t]{7.5cm}{a,b,c,d,e,f,g}{a2=1,b3=b,bc2=b,c3=c,bd=bc,cd=b2,d3=d,be=bc,ce=b2,e2=de,bf=ab,cf=ab2c,d2f=f,f2=b2,b2g=g,c2g=g,dg=cg,eg=cg,fg=ag,g2=b2}
\vspace{0.25cm} \\
\P = & \plist[t]{10cm}{a,b2,ac,ac2,d,ad2,e,ade,adf}
\end{array}\]

\xpretend{12}{a b ab a c ab b ab2 d b bc e ab2 b abc ab2 d2e ab b ade b2c bc abc b2c f b g ab2c b2c abc b ab2 g bc abc b2c ab2 b ab ab2 b2c abc b ab2 g b abc b2c ab2 b g ab2 b2c abc b ab2 b2c b abc b2c ab2 b g ab2 b2c abc b ab2 g bc abc b2c ab2 b g ab2 b2c abc b ab2 g b abc b2c ab2 b g ab2 b2c abc b ab2 g b abc b2c ab2 b g ab2 b2c abc b ab2 g b abc b2c ab2 b g ab2 b2c abc b ab2 g b abc b2c}

\caption{Quotient presentation and pretending function for mis\`ere \textsc{Kayles} to heap 120.}
\label{figure:kayles}
\end{figure}

Now examine the $\Phi$-values $\Phi(H_n) \in \Q$.  We observe that
\[
\Phi(H_{n+12})=\Phi(H_{n}),\textrm{ for }71\leqslant
n\leqslant120-12.
\]
This situation is much like the periodicity of $\mathscr{G}$-values that we observed in normal play.

The following notation will be very useful; it applies to $\textsc{Kayles}$ as well as to an arbitrary octal game~$\Gamma$.  Denote by:
\begin{itemize}
\item $\A$ the set of all positions, $\A=\cl(H_{0},H_{1},\ldots)$.
\item $\Q(\Gamma) = \Q(\A)$.
\item $\A_{n}$ the set of all positions with no heap larger than~$n$, $\A_{n}=\cl(H_{0},\ldots,H_{n})$.
\item $\Q_{n}(\Gamma)=\Q(\A_{n})$, the n$^{\textrm{th}}$ \emph{partial quotient} for~$\Gamma$.
\end{itemize}
We've computed~$\Q_{120}(\textsc{Kayles})$, and
the quotient map $\Phi_{120}:\A_{120}\longrightarrow \Q_{120}$, and
found that it is periodic past a certain point.

\subsection*{A Brief Digression}

In a moment we will state a mis\`ere version of the periodicity theorem.  We first pause to consider some potential difficulties.

\begin{remark}
Suppose we computed~$\Q_{n}$. Now we throw $H_{n+1}$ into the
quotient. There might be games $G,K\in \A_{n}$ such
that~$G\equiv_{\A_{n}}K$ but are distinguished by~$H_{n+1}$. When
this happens, we have~$\Phi_{n}(G)=\Phi_{n}(K)$,
but~$\Phi_{n+1}(G)\neq\Phi_{n+1}(K)$.
\end{remark}

\vspace{0.15cm}

This remark shows that we must be careful not to confuse the partial quotients of $\Gamma$ with its full quotient.

Note that in normal play, there is no such concern.  Given a set of games~$\A$, it is possible to define \emph{normal equivalence modulo~$\A$} in exactly the same way we've defined mis\`ere equivalence modulo~$\A$.  However, in normal play it will always be the case that $G \equiv_\A K$ \emph{if and only if} $G = K$.  That is, in normal play, local and global equivalence coincide.  (To see this, observe that if $G \neq K$ in normal play, then $G$ and $K$ must have different Grundy values, so $G+G$ and $G+K$ have different outcomes.  So if $G$ and $K$ are distinguished by anything, then they must be distinguished locally, by $G$ itself.)  So, although the sorts of localizations we're discussing are perfectly applicable to normal play, they don't provide any further resolution (and in a sense, they don't need to, because normal play is simple enough to begin with).

Let us consider another difference between normal play and mis\`ere play.  Consider a finitely generated set~$\A$. In \emph{normal play}, there
can be only finitely many~$\mathscr{G}$-values represented. To see this, let
$H_1,\ldots,H_n$ generate $\A$.  Then the~$\mathscr{G}$-values
represented by~$\A$ are bitwise exclusive-or's
of~$\mathscr{G}(H_1),\dots,\mathscr{G}(H_n)$, but these are bounded.

What about mis\`{e}re play? Is~$\Q(\A)$ finite?
Answer: Not in general.  Later in this course we will see an example of an infinite, finitely generated quotient.  Our picture of such quotients is still very hazy.  In fact, the following question is still open.

\begin{problem}
Specify an algorithm to determine whether~$\Q(\A)$ is infinite,
assuming~$\A$ is finitely generated.
\end{problem}

We'll say more about this later in the course.  Finally, now is as good a time as any to interject the following remark:

\begin{remark}
All monoids we consider in this course are commutative.  Sometimes I will slip and say ``monoid'' when I really mean ``commutative monoid.''
\end{remark}

\subsection*{Periodicity}

We now return to the setting of an octal game $\Gamma$ with heaps $H_n$.

\vspace{0.15cm}\noindent
\textbf{Recall:}
Periodicity theorem for normal play:
\begin{quote}
Let~$\Gamma$ be an octal game with last non-zero code digit~$k$.
Suppose there are integers~$n_0,p$ such
that~$\mathscr{G}(H_{n+p})=\mathscr{G}(H_{n})$
for~$n_0\leq{}n<2n_0+p+k$.  Then in fact
\[
\mathscr{G}(H_{n+p})=\mathscr{G}(H_{n})\textrm{ for all }n\geq{}n_0.
\]
\end{quote}

\begin{theorem}[Periodicity Theorem for Mis\`{e}re Play]
Let~$\Gamma$ be an octal game with last non-zero code digit~$k$.
Fix~$n_0,p$ and let~$M=2n_0+2p+k$. Let~$(\Q_M,\P_M)=\Q_M(\Gamma)$.
Suppose that~$\Phi_M:A_M\longrightarrow{}\Q_M$, and
that~$\Phi_M(H_{n+p})=\Phi_M(H_n)$ for~$n_0\leq{}n<2n_0+p+k$.  Then
in fact
\[
\Q(\Gamma)\cong{}\Q_M(\Gamma),
\]
and
\[
\Phi(H_{n+p})=\Phi(H_n)\textrm{ for all }n\geq{}n_0.
\]
\end{theorem}

\begin{proof}
Recall the proof for normal play. By induction on~$n$:

\bigskip

\bigskip

\begin{tabular}[c]{|c|c|c|r}
\cline{1-3}
\makebox[3cm]{}&$\times$&\makebox[4cm]{}&\makebox[3.7cm]{$H_n$}\\
\cline{1-3}
\end{tabular}

\bigskip

\begin{tabular}[c]{|c|c|c|r}
\cline{1-3}
\makebox[3cm]{}&$\times$&\makebox[5cm]{}&\makebox[2cm]{$H_{n+p}$}\\
\cline{1-3}
\end{tabular}

\begin{tabular}[t]{c c c r}
$\underbrace{\makebox[3.1cm]{}}$&\makebox[0.06cm]{}&$\underbrace{\makebox[5.3cm]{}}$&\\
$a$&&$b$&\\
\end{tabular}

\bigskip

$H_{n+p}\longrightarrow{}H_a+H_b$ is a typical move from~$H_{n+p}$.
We chose the upper bound of our induction base case to be large
enough that one of~$a,b\geq{}n_0+p$. Assume w.l.o.g that it's~$b$.
But then $\mathscr{G}(H_{b-p})=\mathscr{G}(H_{b})$,
so~$\mathscr{G}(H_a+H_b)=\mathscr{G}(H_a+H_{b-p})$. We conclude that
the options of~$H_n$ and~$H_{n+p}$ represent exactly the
same~$\mathscr{G}$-values. But~$\mathscr{G}$-values are computed by
the mex rule, so this
implies~$\mathscr{G}(H_{n+p})=\mathscr{G}(H_{n})$.

To prove the periodicity theorem for mis\`{e}re play, we can use
exactly the same argument to show that the options of~$H_n,H_{n+p}$
represent exactly the same~$\Phi_M$-values. So the proof now depends
only on the following lemma.

\begin{lemma}
\label{lemma:phiequivalence}
Suppose~$\A$ is a closed set of games, and~$G$ is a game all of whose
options are in~$\A$. Assume that, for some $H \in \A$,
\[\left\{\Phi(G'):G'\textrm{ is an option of
}G\right\}=\left\{\Phi(H'):H'\textrm{ is an option of }H\right\}.\]
Then $\Q(\A\cup\{G\})\cong{}\Q(\A)$ and $\Phi(G)=\Phi(H)$.
\end{lemma}

\emph{Assuming} Lemma~\ref{lemma:phiequivalence}, the proof of the periodicity theorem is complete.  For we can go by induction to show that
\[
Q_M(\Gamma)\cong{}Q_{M+1}(\Gamma)\cong{}Q_{M+2}(\Gamma)\cong\ldots,
\]
and that the resulting $\Phi$-values are periodic.
\end{proof}


\subsection*{Bipartite Monoids}

Although we could prove Lemma~\ref{lemma:phiequivalence} directly, it will be easier after we introduce a suitable abstraction of the mis\`ere quotient construction.  Since the abstract setting is also useful in other situations, this is worth the effort.

\begin{definition}
A \emph{bipartite monoid} is a pair $\QP$ where $\Q$ is a commutative monoid, and $\P \subset \Q$ is some subset. We will usually write b.m. for bipartite monoid.
\end{definition}

\begin{definition}
Let $\QP$ be a b.m. $x,y \in \Q$ are said to be \emph{indistinguishable} if, for all $z \in \Q$,
\[xz \in \P \Longleftrightarrow yz \in \P.\]
\end{definition}

\begin{definition}
A b.m. $\QP$ is \emph{reduced} if the elements of $\Q$ are pairwise distinguishable.  We write r.b.m. for reduced bipartite monoid.
\end{definition}

\begin{examplex}
Every mis\`ere quotient is a r.b.m.
\end{examplex}

\begin{proof}
Suppose $[G]_{\equiv_\A}$ and $[H]_{\equiv_\A}$ are indistinguishable.  Then for any $X \in \A$,
\[[G] + [X] \in \P \Longleftrightarrow [H] + [X] \in \P.\]
Therefore $o^-(G+X) = o^-(H+X)$ for all $X \in \A$, so $[G] = [H]$.
\end{proof}

\begin{examplex}
If $\A$ is a closed set of games, and $\mathscr{B}$ is the set of mis\`ere \PPos{}s of $\A$, then $(\A,\mathscr{B})$ is a bipartite monoid.  The same is true if we take $\mathscr{B}$ to be the set of normal \PPos{}s of $\A$.
\end{examplex}

\begin{definition}
A function $f:\QP \to (\mathcal{S},\mathcal{R})$ is a \emph{bipartite monoid homomorphism} if $f:\Q \to \mathcal{S}$ is a monoid homomorphism,  and for every $x \in \Q$, we have $x\in\P$ iff $f(x)\in \mathcal{R}$.
\end{definition}

\begin{definition}
Let $\QP$ and $(\mathcal{S},\mathcal{R})$ be bipartite monoids.  $(\mathcal{S},\mathcal{R})$ is a \emph{quotient} of $\QP$ iff there is a surjective homomorphism $f:\QP \to (\mathcal{S},\mathcal{R})$.
\end{definition}

\begin{definition}
Let $\QP$ be a b.m.  Define a relation $\rho$ on $\Q$ by  $x \rho y$ iff $x$ and $y$ are indistinguishable.
\end{definition}

\begin{exercise}
Show that $\rho$ is an equivalence relation, and that the equivalence classes modulo $\rho$ form a bipartite monoid.
\end{exercise}

\begin{definition}
The \emph{reduction} of $\QP$ is the bipartite monoid of equivalence classes modulo $\rho$.  We denote it by $(\Q',\P')$.
\end{definition}

\begin{exercise}
Show that $(\Q',\P')$ is reduced and is a quotient of $\QP$.
\end{exercise}

\begin{examplex}
Let $\A$ be a closed set of games, and let $\mathscr{B}$ be the set of mis\`ere \PPos{}s in~$\A$.  Then the mis\`ere quotient $\Q(\A)$ is the reduction of $(\A,\mathscr{B})$.
\end{examplex}

\vspace{0.15cm}
The following proposition is extremely useful.

\begin{proposition}
\label{proposition:uniquereduction}
Suppose $\QP$ is a b.m. with reduction $(\Q',\P')$.  Let $(\mathcal{S},\mathcal{R})$ be any quotient of $\QP$, via $f : \QP \to (\mathcal{S},\mathcal{R})$, and let $(\mathcal{S}',\mathcal{R}')$ be its reduction.  Then there is an isomorphism $i:(\Q',\P') \to (\mathcal{S}',\mathcal{R}')$ making the following diagram commute:
\[\xymatrix{
\Q \ar@{>>}[d] \ar^{f}[r] & \mathcal{S} \ar@{>>}[d] \\
\Q' \ar@{=}^{i}[r] & \mathcal{S}'
}\]
\end{proposition}

\begin{proof}
Let $\rho$ be the reduction relation on $\Q$ ($\Q' = \Q/\rho$), and let $\tau$ be the reduction relation on $\mathcal{S}$ ($\mathcal{S}' = \mathcal{S}/\tau$).

Now for $x,y \in \Q$, we have:
\begin{eqnarray*}
[x]_\rho = [y]_\rho & \textrm{iff} & xz \in \P \Leftrightarrow yz \in \P \textrm{ for all } z \in \Q \\
& \textrm{iff} & f(xz) \in \mathcal{R} \Leftrightarrow f(yz) \in \mathcal{R} \textrm{ for all } z \in \Q \\
& \textrm{iff} & f(x)w \in \mathcal{R} \Leftrightarrow f(y)w \in \mathcal{R} \textrm{ for all } w \in \mathcal{S} 
\qquad \textrm{(since $f$ is surjective)} \\
& \textrm{iff} & [f(x)]_\tau = [f(y)]_\tau
\end{eqnarray*}

So we may define the map $i$ by $i([x]_{\rho}) = [f(x)]_{\tau}$. We just showed that $i$ is well-defined and one-to-one.  Since $f$ is surjective, so is $i$, and it follows that $i$ is an isomorphism.  Commutativity of the diagram follows trivially from the definition of $i$.
\end{proof}

\begin{corollary}
\label{corollary:uniquereduction}
Every bipartite monoid has exactly one reduced quotient (up to isomorphism).
\end{corollary}

Let us see why this is important.  Let $\A$ be a closed set of games, and $\mathscr{B}$ the set of mis\`ere $\mathscr{P}$-positions in $\A$.  Then the mis\`ere quotient $\Q(\A)$ is the reduction of $(\A,\mathscr{B})$.  Therefore, suppose we have some putative quotient $\QP$, and we want to assert that it is $\Q(\A)$.  We just need to show that:
\begin{enumerate}
\item[(a)] $\QP$ is reduced; and
\item[(b)] $\QP$ is a quotient of $(\A,\mathscr{B})$.
\end{enumerate}
By Proposition~\ref{proposition:uniquereduction}, these conditions imply that $\QP \cong \Q(\A)$.  We can therefore avoid the exhaustive analysis used to construct $\mathcal{T}_2$ during the previous lecture.

\newpage
\lecture{4}{November 29, 2006}{More Examples}{Shai Lubliner \& Ohad Manor}

\subsection*{Proof of Lemma~\ref{lemma:phiequivalence}}

We now prove Lemma~\ref{lemma:phiequivalence}, thus completing the proof of the Periodicity Theorem.

\begin{definition}
Suppose $\A$ is a set of games, and $G$ is a game all of whose options are in $\A$.  Define
\[\Phi''G = \{\Phi(G') : G' \textrm{ is an option of } G \}.\]
\end{definition}

(This definition includes the case when $G \in \A$.)

\begin{lemma}
\label{lemma:induction}
Suppose $\A$ is a closed set of games and $\QP$ a r.b.m. The following are equivalent:
\begin{enumerate}
\item[(i)] $\QP \cong \Q(\A)$;
\item[(ii)] There exists a surjective \emph{monoid} homomorphism $\Phi : \A \to \Q$, such that for all $G \in \A$,
\[\Phi(G) \in \P \iff G \neq 0 \textrm{ and } \Phi(G') \not\in \P \textrm{ for every option } G' \textrm{ of } G.\]
\end{enumerate}
\end{lemma}

\begin{proof}
(i) $\Rightarrow$ (ii): Let $\Phi$ be the quotient map $\A \to \Q(\A)$.  We know that, for all $G$,
\[G \textrm{ is a \PPos{} } \iff G \neq 0 \textrm{ and every } G' \textrm{ is an \NPos}.\]
But since $\Phi$ is a homomorphism of bipartite monoids, we have
\[X \textrm{ is a \PPos{} } \iff \Phi(X) \in \P, \textrm{ for all } X \in \A,\]
and the conclusion follows immediately.

\vspace{0.15cm}\noindent
(ii) $\Rightarrow$ (i): By Corollary~\ref{corollary:uniquereduction}, $\Q(\A)$ is the unique reduced quotient of $(\A,\B)$ (where $\B$ is the set of \PPos{}s in $\A$).  Thus it suffices to show that $\Phi$ is a homomorphism of bipartite monoids, since this implies that $\QP$ is a quotient of $(\A,\B)$.  So we must prove the following, for all $G \in \A$:
\[G \textrm{ is a \PPos{} iff } \Phi(G) \in \P.\]
Now by induction on $G$ (i.e., on the height of the game tree of $G$), we may assume that
\[G' \textrm{ is a \PPos{} iff } \Phi(G') \in \P,\]
for all options $G'$ of $G$.  But now:
\[\begin{array}{rcll}
\Phi(G) \in \P & \textrm{iff} & G \neq 0 \textrm{ and } \Phi(G') \not\in \P \textrm{ for all } G' & \textrm{(by assumption)} \\
& \textrm{iff} & G \neq 0 \textrm{ and every } G' \textrm{ is an \NPos} & \textrm{(by induction)} \\
& \textrm{iff} & G \textrm{ is a \PPos} & \textrm{(by definition of \PPos).}
\end{array}\]
This proves the lemma.
\end{proof}

\begin{proof}[Proof of Lemma~\ref{lemma:phiequivalence}]
Assume $\A$ is a closed set of games, all options of $G$ are in $\A$ and $\Phi''G = \Phi''H$ for some $H \in \A$.  We must show that $\Q(\A \cup \{G\}) \cong \Q(\A)$ and $\Phi(G) = \Phi(H)$.

Define $\Phi^{+}: \cl(\A\cup\{G\}) \rightarrow \Q$ by:
\begin{itemize}
\item $\Phi^{+}(G) = \Phi(H)$; 
\item $\Phi^+(Y) = \Phi(Y)$ for all $Y \in \A$.
\end{itemize}
If we regard $G$ as a free generator of the monoid $\cl(\A \cup \{G\})$ over $\A$, then this defines a monoid homomorphism.  So we just need to show that $\Phi^+$ satisfies condition (ii) of Lemma~\ref{lemma:induction}.

Fix $X \in \cl(\A \cup \{G\})$. We can write $X = n \cdot G + Y$ for some $n \geq 0$ and $Y \in \A$.  The case $n = 0$ is already known, so we can assume $n \geq 1$.  Let $W = n \cdot H + Y$; clearly $\Phi^+(X) = \Phi^+(W)$.

Now consider an option $X'$ of $X$.  We have $X' = n \cdot G + Y'$ or $(n-1) \cdot G + G' + Y$.
\begin{itemize}
\item If $X' = n \cdot G + Y'$, then $\Phi^+(X') = \Phi^+(n \cdot H + Y')$, which is an option of $W$.
\item If $X' = (n\Neg 1) \cdot G + G' + Y$, then $\Phi^+(X') = \Phi^+((n\Neg 1) \cdot H + G' + Y)$.  But since $\Phi''G = \Phi''H$, there must be some $H'$ with $\Phi^+(H') = \Phi^+(G')$.  So $\Phi^+(X') = \Phi^+((n\Neg 1) \cdot H + H' + Y)$, again an option of~$W'$.
\end{itemize}
This shows that $(\Phi^+)''X \subset (\Phi^+)''W$, and an identical argument shows that $(\Phi^+)''W \subset (\Phi^+)''X$.  But since $W \in \A$, we know that
\[\Phi(W) \in \P \iff W \neq 0 \textrm{ and } \Phi(W') \not\in \P \textrm{ for all } W'.\]
Since $\Phi^+(X) = \Phi^+(W)$ and $(\Phi^+)''X = (\Phi^+)''W$, we have
\[\Phi^+(X) \in \P \iff W \neq 0 \textrm{ and } \Phi^+(X') \not\in \P \textrm{ for all } X'.\]
This satisfies Lemma~\ref{lemma:induction}(ii) except for the condition $W \neq 0$.  But if either of $G,H$ is identically $0$, then both must be, since $\Phi''G = \emptyset$ iff $\Phi''H = \emptyset$.  Therefore $W \neq 0$ iff $X \neq 0$, and we are done.
\end{proof}

\subsection*{Further Examples}

The partial quotients of \textsc{Nim} are fundamental examples, and we denote them by $\mathcal{T}_n$.
\begin{itemize}
\item $\T_0 = \Q(0)$;
\item $\T_1 = \Q(*)$;
\item $\T_2 = \Q(*2)$;
\item $\T_n = \Q(*2^{n-1})$.
\end{itemize}
Here are their presentations:
\begin{itemize}
\item $\T_0 = \{ 1 \}$;  $\P =  \emptyset$
\item $\T_1 = \<a~|~a^2=1\>$; $\P = \{a\}$
\item $\T_2 = \<a,b~|~a^2=1,\ b^3=b\>$; $\P = \{a,b^2\}$
\item $\T_3 = \<a,b,c~|~a^2 = 1,\ b^3 = b,\ c^3 = c,\ b^2 = c^2\>$; $\P = \{a,b^2\}$
\item $\T_n = \<a,b_1,b_2,\ldots,b_{n-1}~|~a^2=1,\ b_i^3=b_i,\ b_1^2=b_2^2=\cdots=b_{n-1}^2\>$; $\P = \{a,b_1^2\}$
\end{itemize}
To find $\Phi(*m)$ (in any of the $\T_n$), write $m$ in binary, as $\cdots \epsilon_3\epsilon_2\epsilon_1\epsilon_0$, and we have
\[\Phi(*m) = a^{\epsilon_0} \cdot b_1^{\epsilon_1} \cdot b_2^{\epsilon_2} \cdot \cdots \cdot b_n^{\epsilon_n}.\]
For example, in $\T_4$, we have
\[\Phi(*4) = b_2,\ \Phi(*5) = ab_2,\ \Phi(*6) = b_1b_2,\ \Phi(*7) = ab_1b_2,\ \Phi(*8) = b_3\]
Notice that we always have
\[b_1^2 = b_2^2 = \cdots = b_{n-1}^2\]
Denote this element by $z$.  $z$ represents the sum $*m + *m$, for \emph{any} \textsc{Nim}-heap with $m \geq 2$.  In fact, it represents any \textsc{Nim} position of $\mathscr{G}$-value $0$, provided it has at least one heap of size $\geq 2$.

\subsection*{The Structure of $\T_n$}

Let's write out the elements of $\T_3$.
\[\T_3 = \{1, a, b_1, ab_1, b_2, ab_2, b_1 b_2, ab_1 b_2, z, az\}\]
Consider the subset
\[\mathcal{K} = \{ b_1, ab_1, b_2, ab_2, b_1 b_2, ab_1 b_2, z, az\}\]
Observe that $z \cdot z = z$, $z \cdot b_1 = b_1$, and $z \cdot b_2 = b_2$.  Therefore $z$ is an identity of $\mathcal{K}$ and $x^2 = z$ for all $x \in \mathcal{K}$.  So $\mathcal{K}$ is a group, and we have
\[\mathcal{K} \cong \mathbb{Z}_2^3.\]
In fact $\mathcal{K}$ behaves just like normal play $\mathscr{G}$-values: it has eight elements, corresponding one-to-one with \textsc{Nim} positions of $\mathscr{G}$-value $0$ through $7$.

Recall the strategy for mis\`ere \textsc{Nim}:
play exactly like in normal \textsc{Nim}, unless
your move would leave only heaps of size 0 or 1. In that case, play
to leave an odd number of heaps of size 1.

$\mathcal{K}$ corresponds to the ``exactly like normal \textsc{Nim}'' clause of this strategy: it is isomorphic to the normal-play quotient of $*4$.  The two elements $1$ and $a$ correspond to the ``unless'': they represent positions with all heaps of size $\leq 1$.

Note that every $\T_n$, for $n \geq 2$, can be written as $\mathcal{K} \cup \{1,a\}$, where $\mathcal{K} \cong \mathbb{Z}_2^n$.  $\mathcal{K}$ is called the \emph{kernel} of the monoid, and in the next lecture we will see how to generalize it.

In particular we have:
\begin{itemize}
\item $|\T_0| = 1$
\item $|\T_1| = 2$
\item $|\T_n| = 2^n+2$ for all $n \geq 2$
\end{itemize}

We can also define the full quotient of \textsc{Nim}:
\[\T_\infty = \Q(0,*,*2,*3,*4,\ldots) \cong \<a, b_1, b_2~|~a^2 = 1,\ b^3_i = b_i,\ b^2_1 = b^2_2 = \ldots\> \qquad \P = \{a,b_1^2\}\]
Remember that normal-play $\mathscr{G}$-values look like
\[\bigoplus_{\mathbb{N}} \mathbb{Z}_2\]
Well, we can write $\T_\infty = \mathcal{K}_\infty \cup \{1,a\}$ in exactly the same way, and we have $\mathcal{K}_\infty \cong \bigoplus_{\mathbb{N}} \mathbb{Z}_2$.

\subsection*{Tame and Wild Quotients}

\begin{definition}
\label{def:tameness}
A set $\A$ is \emph{tame} iff $\Q(\A) \cong \T_n$ for some $n \in \mathbb{N} \cup \{\infty\}$.  Otherwise it is \emph{wild}.
\end{definition}


Not all quotients are tame:

\begin{examplex}
Let $G = \Q(\Star 2_\sh 320)$, where $\Star 2_\sh320 = \{0,*2,*3,*2_\sh\}$ and $*2_\sh = \{*2\}$.  We have
\[\Q(G) \cong \<a, b, t~|~a^2 = 1,\ b^3 = b,\ t^2 = b^2,\ bt = b\>;\ \P = \{a,b^2\}\]
This quotient is called $\mathcal{R}_8$.  It is very common; many octal games have quotient $\mathcal{R}_8$, including (for example) 0.75.  In fact, it can be shown that $\mathcal{R}_8$ is the smallest quotient except for $\T_0$, $\T_1$, $\T_2$.  The quotient map is given by (writing $z = b^2$, as before)
\[\Phi(*) = a,\ \Phi(*2) = b,\ \Phi(*3) = ab,\ \Phi(*2_\sh) = z,\ \Phi(G) = at.\]
\end{examplex}

Notice that $\mathcal{R}_8$ is just $\T_2$ with two extra elements:
\[\Q = \{ \underbrace{1, a, \underbrace{b, ab, z, az}_{\mathcal{K}}}_{\T_2}, t, at\}\]
Now $\mathcal{K} \cong \mathbb{Z}_2^2$, and $\{1,a\}$ is a (separate) isomorphic copy of $\mathbb{Z}_2$.  But $\{t,at\}$ is not a group, because $t^2 = z \in \mathcal{K}$.

The right picture of $\mathcal{R}_8$ is this: it is the union
\[\mathcal{K} \cup \{1,a\} \cup \{t,at\},\]
where $\mathcal{K}$ and $\{1,a\}$ are two disjoint groups, and $\{t,at\}$ are two extra elements that are ``associated'' with $\mathcal{K}$.  We'll say more about this in the next lecture.

\subsection*{General Structure}

\begin{lemma}
\label{lem:heredity}
Suppose that $\A$ is hereditarily closed, $\A \neq \emptyset$, and $\A \neq \{0\}$.  Then necessarily $* \in \A$. 
\end{lemma}

\begin{proof}
\label{prf:heredity}
$*$ is the only game whose only option is $0$.
\end{proof}

\begin{proposition}
\label{pro:position}
Let $\QP$ be any nontrivial mis\`ere quotient.  Then  for all $x \in \Q$, there is some $y \in \Q$ with $xy \in \P$.
\end{proposition}

\begin{proof}
\label{prf:position}
Write $\QP = \Q(\A)$ and choose $ G \in \A $ with $ \Phi(G) = x$.  First suppose $G = 0$.  Then $x = 1$.  By the assumption of nontriviality, we have $\A \neq \{0\}$, so by the previous lemma $* \in \A$.  But $\Phi(*) \in \P$ and $1 \not\in \P$, so we can take $y = \Phi(*)$.

Now assume $G \neq 0$, and consider $G + G$.  If it is a \PPos, then we are done, with $y = x$.  Otherwise, some option of $G + G$ must be a \PPos, say $G + G'$.  So we can take $y = \Phi(G')$.
\end{proof}

\begin{proposition}
\label{pro:option}
For any $G$ and any option $G^{\prime}$, $\Phi(G) \neq \Phi(G^{\prime})$.
\end{proposition}

\begin{proof}
\label{prf:option}
Exercise.  (Hint: Use the previous proposition.)
\end{proof}

\begin{proposition}
\label{pro:star}
If $\A$ is nontrivial and $G \in \A$, then $G \not\equiv_{\A} G+*$ 
\end{proposition}

\begin{proof}
\label{prf:star}
By Proposition~\ref{pro:position}, there is a game $H \in A$ such that $G+H$ is a \PPos.  But then $G + H + *$ is an \NPos, so $H$ distinguishes $G$ from $G+*$.
\end{proof}

\begin{corollary}
\label{cor:even}
Every nontrivial mis\`ere quotient has even order.
\end{corollary}

\begin{proof}
\label{prf:even}
Exercise.  (Hint: Consider the mapping $x \mapsto ax$.)
\end{proof}

In fact, one can prove the following facts.
\begin{itemize}
\item $\T_1$ is the only quotient of order 2.  (Immediate from Lemma~\ref{lem:heredity})
\item There are no quotients of order 4.  (Proved in~\cite{siegel_200Xd})
\item $\T_2$ is the only quotient of order 6.  (Also proved in~\cite{siegel_200Xd})
\item $\mathcal{R}_8$ is the only quotient of order 8.  (Much harder to prove; see~\cite{siegel_200Xe})
\end{itemize}

\newpage
\lecture{5}{November 30, 2006}{Further Topics}{Shiri Chechik \& Menachem Rosenfeld}

\noindent In this lecture, we will discuss four interesting problems, most of which have not yet been solved completely. We will also discuss the structure of finite commutative monoids.

\section*{Four interesting problems}

\subsection*{1. Infinite Quotients}

We can think of infinite quotients as belonging to either one of two categories: Those that are finitely generated, and those that are not. We have already seen one infinite quotient, $\mathcal{T}_{\infty} = \mathcal{Q}(0, *, *2, \ldots)$. It is not finitely generated. Every one of its finitely generated submonoids is finite, and it is built up from these finite quotients. It is therefore not an interesting quotient to study.

\begin{figure}[hb]
\[
\begin{array}{cc}
\xymatrix{&E\ar[dl]\ar[dr]\\0&&D\ar[dl]\ar[dr]\\&0&&C\ar[d]\\&&&\ast2}
&
\xymatrix{\framebox{A}\ar@2{->}[d]&\framebox{B}\ar@2{->}[l]\ar@2{->}[d]&\framebox{C}\ar@2{->}[l]&\framebox{D}\ar@2{->}[d]\ar@2{->}[l]&\framebox{E}\ar@2{->}[d]\ar@2{->}[l]\\&&&&}
\\
(a)&(b)
\end{array}
\]
	\caption{Two representations for the game $E=*(2_\#0)0$. ($a$) The game tree of $E$. ($b$) A visual representation of $\mathrm{cl}(E)$.}
	\label{fig:ETree}
\end{figure}
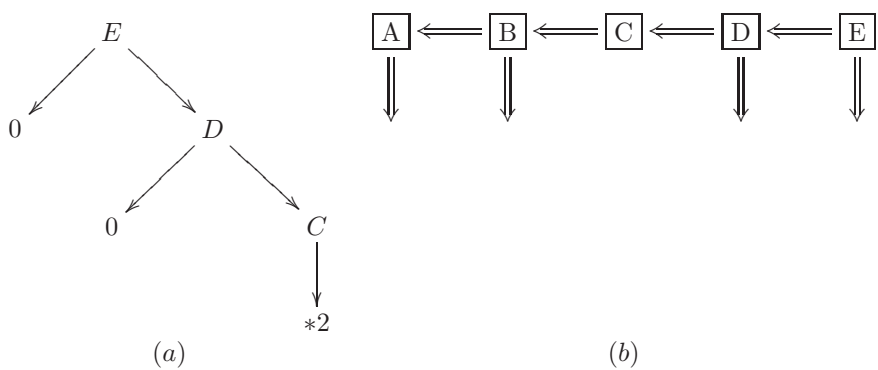

There also exist finitely generated infinite quotients. We can find an example of this by denoting
\[A = *,\ B = *2,\ C = \{B\} = *2_\sh,\ D = *2_\sh0 = \{C, 0\}, \textrm{ and } E = \{D, 0\} = *(2_\sh0)0.\]
Figure  \ref{fig:ETree}($a$) shows the game tree of $E$.

Denoting $\mathscr{A} = \textrm{cl}(E)$, a visual way to understand a game in $\mathscr{A}$ is suggested in Figure \ref{fig:ETree}($b$); for every game, there are several coins in every box, and a move consists of moving a coin along an arrow (either one step to the left, or from boxes other than $C$, outside the game board). The last player to move loses.

\begin{figure}
\centering
\begin{tabular}{@{}c@{\hspace{-8pt}}c@{\hspace{-8pt}}c@{\hspace{-8pt}}c@{}}
& \hspace{8pt} $k=0$ & \hspace{8pt} $k=1$ & \hspace{8pt} $k=2$ \\
\input ppostable.tex
\end{tabular}

\[\bullet = \mathscr{P}; \qquad \textrm{x} = \mathscr{P} \textrm{ iff } j \geq 2; \qquad \circ = \mathscr{P} \textrm{ iff } j < 2\]

\vspace{-11pt}

\caption{Schematic of the \PPos{}s for $\cl(\Star(2_\# 0)0)$ with $k \leq 2$.}
\label{figure:infquotient}
\end{figure}

As it turns out, $|\Q(E)| = \infty$, but every game with a smaller tree has a finite quotient.  So $E$ is in some sense the \emph{simplest} game that gives rise to an infinite quotient.  To understand why the quotient is infinite, first note that every $X \in \A$ can be written as $X = iA+jB+kC+lD+mE$.  In~\cite[Section 6]{siegel_200Xd}, we compute the outcome of every such $X$.  It turns out that when $k \geq 3$, the outcomes follow a simple rule: $o^-(X)=\mathscr{P} \iff i+l$ and $j+m$ are both even.  However, when $k \leq 2$, the outcomes can be quite erratic.  See Figure~\ref{figure:infquotient}.  Each table represents the outcomes for a particular choice of $(i,j,k)$.  Within each table, there is a dot at (row $m$, column $l$) iff $iA+jB+kC+lD+mE$ is a \PPos.

Inspecting this figure, we can see that the structure of the \PPos{}s is very complicated. For example, for $i=j=k=l=0$, $X=mE$ is a \PPos{} $\iff$ $m \in \{1,4,7,10,12,14,16,\ldots\}$.

To see that the quotient is infinite, consider the case $i=j=k=0$.  For sufficiently large odd $l$, we have that $lD+mE$ is a \PPos{} iff $m=l+7$. This means that the $lD$'s are pairwise distinguishable. 

It was mentioned in a previous lecture that infinite quotients are still poorly understood.  We still cannot solve the following problem.

\begin{problem}
Specify an algorithm to determine whether a quotient is infinite.
\end{problem}

Of course, we'd really like to know much more about $\Q(\A)$ than merely \emph{whether} it's infinite.  An old theorem about commutative semigroups guarantees that this is possible:

\begin{theorem}[R\'edei]
Every finitely generated commutative semigroup is finitely presented.
\end{theorem}

We won't prove R\'edei's Theorem in this course; see~\cite{grillet_2001,redei_1965}.  It makes the following question meaningful.

\begin{problem}
Specify an algorithm to compute the \emph{presentation} of~$\Q(\A)$ (even if~$\Q$ is infinite), assuming~$\A$ is finitely generated.
\end{problem}

In particular, the following would be a good start.

\begin{problem}
Give a presentation for $\mathcal{Q}(E)$.
\end{problem}

Note: when we proved the Periodicity Theorem, at no point did we assume that the partial quotients are finite.  Thus the Periodicity Theorem applies perfectly well to octal games whose partial quotients are infinite.  If we could produce an algorithm for computing infinite quotients, then we could (in theory) use the Periodicity Theorem to provide solutions to games with infinite partial quotients.

\subsection*{2. Algebraic Periodicity of Octal Games}
Let $\Gamma$ be an octal game. Then $\mathcal{Q}(\Gamma)$ is uniquely determined by its sequence of partial quotients, \[\left\langle \Q_n(\Gamma):n \in \mathbb{N} \right\rangle.\]
We can ask, when is it determined by only finitely many of these partial quotients?

The periodicity theorem is a good start in trying to answer this question---it happens, for instance, when the sequence stabilizes and we have periodicity.

There are intriguing cases in which the sequence does not stabilize but exhibits a strong regularity, which is called \emph{algebraic periodicity}. This phenomenon is not yet understood well enough for a precise definition to be given. The term is derived from \emph{arithmetic periodicity} in normal play, which means that the sequence is periodic but on each period we add a ``saltus''. For example, if the period is 5 and the saltus is 4, a possible sequence is \[0,4,5,3,2,\ 4,8,9,7,6,\ 8,12,13,11,10,\ \ldots\]

\begin{theorem}
No finite octal game (that is, one with finitely many non-zero digits) can be arithmetic periodic (with non-trivial saltus) in normal play.
\end{theorem}

(Remark: \textsc{Nim} is a trivial example of a non-finite octal game which is arithmetic periodic.)

However, algebraic periodicity is manifested in finite octal games with mis\`ere play. Page 38 of \cite{siegel_200Xd} presents several examples.

Several two-digit octal games for which the normal solution is known, have not yet been solved for mis\`ere play. Of these, 0.54 is the only one for which the solution seems to be in reach---because it appears to be algebraic periodic, which suggests a solution for it.

\begin{problem}
Prove this solution for 0.54.
\end{problem}

\begin{problem}
Formulate a suitable general definition of ``algebraic periodicity'' and prove a theorem that states: If $\Gamma$ is algebraic periodic for sufficiently long, then it continues this period, and we can compute $\mathcal{Q}(\Gamma)$.
\end{problem}

Presumably, this would immediately provide a solution for 0.54, and probably six or eight three-digit octals as well.

\subsection*{3. Generalizations of the Mex Rule}
Suppose we have a quotient map $\Phi:\mathscr{A} \to \mathcal{Q}$. Let $G$ be a game all of whose options are in $\mathscr{A}$. Can we determine, based only on $\Phi''G$, whether $\mathcal{Q}(\mathscr{A} \cup \{G\}) \cong \mathcal{Q}(\mathscr{A})$? If they are isomorphic, can we determine $\Phi(G)$?  (Recall that $\Phi''G$ is defined as $\{\Phi(G') : G'$ is an option of $G\}$.)

By asking these questions, we are essentially looking for a way to generalize the mex rule, which solves them for normal play.

The answer to both question is: Yes! However, more information is needed than what it contained in $\mathcal{Q}$.

Recall that in the previous lecture, we proved a lemma that answers this question in case there is some $H \in \mathscr{A}$ such that $\Phi''G=\Phi''H$. It turns out that we can get a much stronger result. However, this result is beyond the scope of this lecture; see \cite[Section 7]{siegel_200Xd}.

\subsection*{4. Classification}
How many Mis\`ere quotients are there of order $n$ (up to isomorphism)? Table \ref{tab:NumQuotients} displays some of what is known so far.  The results for $n = 14$ and $16$ are tentative.

\begin{table}[h]\centering
	\begin{tabular}{c|ccccccccc}
	Order (n)&2&4&6&8&10&12&14&16&$\cdots$\\\hline\# of quotients&1&0&1&1&1&6&9&50&$\cdots$
	\end{tabular}
	\caption{Number of different quotients for every order}
	\label{tab:NumQuotients}
\end{table}

\emph{Related question:} Can we identify other interesting classification results?  Here is one such result.

It is possible to define the ``tame extension'' $\mathcal{T}(\mathcal{Q},\mathcal{P})$ of an arbitrary quotient $(\mathcal{Q},\mathcal{P})$.  See~\cite{siegel_200Xe} for a precise definition.  It turns out that
\[(\mathcal{Q},\mathcal{P}) \subsetneqq \mathcal{T}(\mathcal{Q},\mathcal{P})\]
but $\mathcal{T}(\mathcal{Q},\mathcal{P})$ adds no new $\mathscr{P}$-position types.  Furthermore,
\[\mathcal{T}_{n+1} = \mathcal{T}(\mathcal{T}_n).\]
We therefore have two families of quotients,
\[\mathcal{T}_2, \mathcal{T}_3, \ldots \mathcal{T}_{\infty}\]
and
\[\mathcal{R}_8, \mathcal{T}(\mathcal{R}_8), \mathcal{T}(\mathcal{T}(\mathcal{R}_8)), \ldots \mathcal{T}^\infty(\mathcal{R}_8),\]
all of which have $|\mathcal{P}| = 2$.  The following result is proved in \cite{siegel_200Xe}.


\begin{theorem}
Every quotient with $|\mathcal{P}| = 2$ is isomorphic to a quotient in one of these two families.
\end{theorem}

So we have:

\begin{center}
\begin{tabular}{ll}
$0, \mathcal{T}(0), \mathcal{T}(\mathcal{T}(0)), \ldots$ & Normal play\\
$\mathcal{T}_2, \mathcal{T}(\mathcal{T}_2), \ldots, \mathcal{T}^{\infty}(\mathcal{T}_2)$ & Tame Mis\`ere play\\
$\mathcal{R}_8, \mathcal{T}(\mathcal{R}_8), \mathcal{T}^2(\mathcal{R}_8), \ldots, \mathcal{T}^{\infty}(\mathcal{R}_8)$ & ``Almost tame'' Mis\`ere play
\end{tabular}
\end{center}

Can we say anything else along these lines?

\section*{The Structure of Finite Commutative Monoids}

Let $\mathcal{Q}$ be any finite commutative monoid, and let $x,y \in \mathcal{Q}$.

\begin{definition}
$x$ \emph{divides} $y$ if $xz=y$ for some $z \in \mathcal{Q}$. In this case, we write $x|y$.
\end{definition}

\begin{definition}
$x$ and $y$ are \emph{mutually divisble} (shorthand : m.d.) if $x|y$ and $y|x$.
\end{definition}

\begin{examplex}
$\mathcal{T}_2 = \left\langle a,b~|~a^2=1,\ b^3=b \right\rangle = \{\underbrace{1,a}_{\textrm{m.d.}},\underbrace{b,ab,z,az}_{\textrm{m.d.}}\}$
\end{examplex}

\begin{exercise}
Show that m.d. is an equivalence relation.
\end{exercise}

\begin{definition}
The \emph{mutual divisibility classes} of $\mathcal{Q}$ are the equivalence classes of $\mathcal{Q}$ under the relation m.d.
\end{definition}

\begin{examplex}
The m.d. classes of $\mathcal{R}_8 = \{1,a,b,ab,z,az,t,at\}$ are $\{1,a\}$, $\{b,ab,z,az\}$ and $\{t,at\}$.
\end{examplex}

\begin{definition}
An element $x\in \mathcal{Q}$ is an \emph{idempotent} if $x^2=x$.
\end{definition}

\begin{examplex}
In $\mathcal{T}_2$ (and also $\mathcal{R}_8$) $1$ and $z$ are the only idempotents.
\end{examplex}

\begin{exercise}
The m.d. class of an idempotent $x$ is a group with $x$ for an identity.
\end{exercise}

\begin{exercise}
If $S$ is a maximal subgroup of $\mathcal{Q}$ (that is, a group which is not contained in any larger subgroup of $\mathcal{Q}$) then it is the m.d. class of its idempotent.
\end{exercise}

Let $z_1, z_2, \ldots, z_k$ be the idempotents of $\mathcal{Q}$ (since $\mathcal{Q}$ is finite, we can enumerate them all). We write $z = z_1z_2\cdots z_k$. We then have $z^2 = z_1^2z_2^2\cdots z_k^2 = z_1z_2\cdots z_k = z$ and $zz_i = z$.

\begin{definition}
The \emph{kernel} of $\mathcal{Q}$ is the m.d.\ class of $z$, and is denoted $\mathcal{K}$.
\end{definition}

We will soon prove the following theorem.

\begin{theorem}
The map $x \mapsto zx$ is a surjective homomorphism from $\mathcal{Q}$ onto $\mathcal{K}$.
\end{theorem}

The kernel $\mathcal{K}$ can be characterized in two ways:
\begin{enumerate}
\item It is the unique group such that there is a surjective homomorphism $f:\mathcal{Q} \to \mathcal{K}$ with the property: If $g:\mathcal{Q} \to \mathcal{D}$ is a homomorphism onto a group $\mathcal{D}$, then there exists an $h:\mathcal{K} \to \mathcal{D}$ which makes the following diagram commute:

\[\xymatrix{\mathcal{Q}\ar[r]^f\ar[dr]_g&\mathcal{K}\ar[d]^h\\&\mathcal{D}}\]

In other words, any map from $\mathcal{Q}$ onto a group $\mathcal{D}$ factors through $f$.
\item$\mathcal{K}$ is the \emph{group of fractions} of $\mathcal{Q}$, that is, it is the group obtained by adjoining formal inverses to $\mathcal{Q}$.
\end{enumerate}

\begin{lemma}
If $y \in \mathcal{Q}$ then for some $r > 0$, $y^r$ is an idempotent.
\end{lemma}
(Note: This does not hold for infinite monoids!)

\begin{proof} Consider the sequence $y,y^2,y^3,y^4,\ldots$. Since $\mathcal{Q}$ is finite, there must be some $n>0$ and some $k>0$ such that $y^n = y^{n+k}$. We then have for every $t \ge 0$, $y^{n+tk} = y^n$. Let $r$ be the unique integer such that $n \le r < n+k$ and $r \equiv 0 \pmod{k}$. Then:
\[y^{2r} = y^{r+tk} = y^{n+tk}y^{r-n} = y^ny^{r-n} = y^r\]
So $y^r$ is an idempotent.
\end{proof}

Note that this idempotent is uniquely determined for any given $y$. Therefore, for any $y$, there is a unique idempotent $x$ such that $y^n=x$ for some $n > 0$. This motivates the following definition:

\begin{definition}
For any idempotent $x \in \mathcal{Q}$, the \emph{Archimedean component} of $x$ is $\{y\in \mathcal{Q} : \exists n(y^n=x)\}$.
\end{definition}

What we have actually shown is that every $y \in \mathcal{Q}$ is a member of a unique archimedean component. Therefore, $\mathcal{Q}$ is partitioned into several Archimedean components. For example, $\mathcal{R}_8$ is partitioned into $\{1,a\}$ and $\{b,ab,z,az,t,at\}$.

We complete the picture by defining a natural partial order on idempotents.

\begin{definition}
For idempotents $x,y \in \mathcal{Q}$, $x \le y \iff xy=x$.
\end{definition}

\begin{examplex}
For any idempotent $x$, $z \le x \le 1$.
\end{examplex}

\begin{theorem}
The idempotents of $\mathcal{Q}$ form a lattice with respect to the relation $\le$.
\end{theorem}

\begin{exercise}
Prove this theorem.  \emph{Hint}: Define
\[x \wedge y = xy\]
and
\[x \vee y = \prod\{w \in \mathcal{Q} : w \textrm{ is an idempotent and } w \ge x,y\}.\]
\end{exercise}

\begin{examples}
In these examples, $[a^b]$ denotes an Archimedean component with $a$ elements contained in the m.d.\ class of the idempotent, and $b$ additional elements.

\begin{center}
\begin{tabular}{c@{\hspace{2cm}}c}
$\mathcal{R}_8$:& \textsc{Kayles}: \vspace{0.25cm}\\
\xymatrix{1~[2]\ar@{-}[d]\\z~[4^2]}&
\xymatrix{
&1~[2]\ar@{-}[dl]\ar@{-}[ddr]\\
[4]~d^2 \ar@{-}[dd]\\
&&c^2~[4]\ar@{-}[ddl]\\
[4^2]~de\ar@{-}[dr]\\
&b^2~[16^8]}
\end{tabular}
\end{center}
\end{examples}

\begin{theorem}
The map $x \mapsto zx$ is a surjective homomorphism $\mathcal{Q} \to \mathcal{K}$.
\end{theorem}

\begin{proof}
We must show that $zx \in \mathcal{K}$ for all $x$; it follows easily that $x \mapsto zx$ is a surjective homomorphism.

Clearly $z|zx$, so we must show that $zx|z$. Let us take an $n > 0$ such that $x^n$ is an idempotent. Then $zx^n=z$ by the definition of $z$, so $(zx)x^{n-1}=z$.
\end{proof}

\begin{corollary}
If $x \in \mathcal{K}$, then $\forall y \in \mathcal{Q}$, $xy \in \mathcal{K}$ (because $x=zx$, so $xy=zxy$).
\end{corollary}

\begin{corollary}
$\mathcal{K} \cap \mathcal{P} \neq \emptyset$ (because by a previous lemma, $\exists x \in \mathcal{Q}$ such that $xz \in \P$, but $xz \in \mathcal{K}$ so $xz \in \mathcal{K} \cap \P$).
\end{corollary}

\begin{definition}
$(\mathcal{Q},\mathcal{P})$ is \emph{normal} if $\mathcal{K} \cap \mathcal{P} = \{z\}$.
\end{definition}

(Remark: the smallest known example of an abnormal quotient is of size 420.)

\begin{definition}
$(\mathcal{Q},\mathcal{P})$ is \emph{regular} if $|\mathcal{K} \cap \mathcal{P}|=1$.
\end{definition}

(Remark: the smallest known example of an irregular quotient is of size over 3000.)

\begin{definition}
A quotient map $\Phi:\mathscr{A} \to \mathcal{Q}$ is \emph{faithful} if, for all $G,H \in \mathscr{A}$,
\[\Phi(G) = \Phi(H) \Rightarrow \mathscr{G}(G) = \mathscr{G}(H).\]
\end{definition}

\begin{problem}
Is every quotient map faithful?
\end{problem}

\begin{theorem}
If $(\mathcal{Q},\mathcal{P})$ is normal and $\Phi$ is faithful, then for all $G,H \in \mathscr{A}$,
\[z\Phi(G) = z\Phi(H) \iff \mathscr{G}(G) = \mathscr{G}(H).\]
\end{theorem}

There is therefore a one-to-one correspondence between elements of $\mathcal{K}$ and normal-play Grundy values of games in $\mathscr{A}$. Furthermore, we can compute the mex function in the kernel. This gives us the following strategy for playing a mis\`ere octal game $\Gamma$: Play as if you were playing normal $\Gamma$, unless your move would take you outside of $\mathcal{K}$. Then pay attention to the fine structure of the mis\`ere quotient.

\begin{examplex}
The octal game 0.414 has not yet been solved for normal play. Nevertheless, we can prove that $\Phi(H_n) \in \mathcal{K}$ for $n>18$, and we can prove that its quotient is one of
\[\mathcal{Q}_{18}, \mathcal{T}(\mathcal{Q}_{18}), \mathcal{T}(\mathcal{T}(\mathcal{Q}_{18})), \ldots,\]
though we do not know which.  The strategy for mis\`ere 0.414 is: play as if you were playing normal 0.414, unless your move would leave only heaps of size $\leq 18$.  Then pay attention to the fine structure of the mis\`ere quotient.
\end{examplex}

\vspace{0.15cm}
One last open problem:

\begin{problem}
Let $\mathcal{S}$ be an arbitrary maximal subgroup of $\mathcal{Q}$. Must $\mathcal{S} \cap \mathcal{P}$ be non-empty?
\end{problem}

\newpage

\section*{Further Reading}

\emph{Mis\`ere quotients for impartial games} \cite{siegel_200Xd}, by Plambeck and Siegel, includes most of the material presented in these notes, and a great deal else as well.  It is the best resource both for additional examples of mis\`ere quotients and for a deeper look at the structure theory.  Plambeck's original paper introducing mis\`ere quotients \cite{plambeck_2005} includes a proof of the periodicity theorem that is somewhat different from the one presented here.  His survey paper \cite{plambeck_200X} provides a nice informal summary of much that is known about mis\`ere games.  The forthcoming paper \cite{siegel_200Xe} dives much more deeply into the structure of mis\`ere quotients.

The most current source of information is the \emph{Mis\'ere Games} website~\cite{miseregames_www}, which includes Plambeck's mis\'ere games blog.  See also~\cite{plambeck_www}.

The \emph{canonical} theory is virtually useless in practice, but nonetheless absolutely fascinating.  It is (essentially) the ``quotient'' obtained by taking $\A$ to be the universe of all mis\`ere games.  See \cite{conway_1976} and the forthcoming paper \cite{siegel_200Xg}.

Finally, perhaps the best way to get acquainted with mis\`ere quotients is to download a copy of \emph{MisereSolver}~\cite{miseresolver} and start experimenting.  It can easily reproduce all the examples in this paper, and of course many more as well.

\newpage

\bibliography{games}

\end{document}